\documentclass[brochure,12pt]{bourbaki}
\usepackage[matrix,arrow]{xy}
\usepackage{amssymb,amsfonts,amsmath,footnote}
\usepackage[applemac]{inputenc}
\usepackage[francais]{babel}
\begin{abstract}
Les groupes sp\'eciaux lin\'eaires entiers, les groupes modulaires de surfaces
et les groupes des automorphismes ext\'erieurs de groupes libres
apparaissent dans de nombreux domaines. Leurs analogies, soulign\'ees en
particulier par les travaux de K.~Vogtmann, font couler beaucoup d'encre.
Dans ce rapport, nous nous concentrerons sur les espaces contractiles sur
lesquels ces groupes agissent de mani\`ere analogue, sur les propri\'et\'es
communes de leurs sous-groupes et sur les propri\'et\'es semblables (ou
envisag\'ees semblables) de leur g\'eom\'etrie asymptotique.
\end{abstract}
\altkeywords{Automorphism group, free group, surface group, special linear group,
group action on trees, Teichm\"uller space, Culler-Vogtmann's Outer space,
asymptotic geometry of groups.}
\begin{altabstract}
The special linear groups, the mapping class groups of surfaces, 
the outer autormorphism groups of free groups appear in numerous 
domains. Their analogies, developped in particular in K.~Vogtmann's 
work, have been written about a lot. In this report, we concentrate 
on the contractible spaces on which these groups act in an analogous 
way, on the common properties of their subgroups, and on the similar 
(or conjecturally similar) properties of their asymptotic geometry.
\end{altabstract}

\title
[Automorphismes de groupes libres]
{Sur les automorphismes de groupes libres\\ 
et de groupes de surface}
\alttitle{On the automorphism groups of free
groups and surface groups}

\date{Juin 2010}
\bbkannee{62\`eme ann\'ee, 2009-2010}
\bbknumero{1023}
\author{Fr\'ed\'eric PAULIN}
\address{Ecole Normale Sup\'erieure\\
DMA, UMR 8553 CNRS\\
45 rue d'Ulm\\
F--75230 PARIS C\'edex 05\\
et\\
Universit\'e Paris-Sud 11\\
UMR 8628 CNRS, B\^at.~425\\
D\'epartement de Math\'ematiques\\
F--91405 ORSAY C\'edex}
\email{frederic.paulin@math.u-psud.fr}

\usepackage[T1]{fontenc}
\usepackage{comment}
\usepackage[dvips]{color}
\usepackage{mathrsfs}

\newcommand{\bdefi}{\begin{defi}}
\newcommand{\edefi}{\end{defi}}
\newcommand{\bprop}{\begin{prop}}
\newcommand{\eprop}{\end{prop}}
\newcommand{\btheo}{\begin{theo}}
\newcommand{\etheo}{\end{theo}}
\newcommand{\blemm}{\begin{lemm}}
\newcommand{\brema}{\begin{rema}}
\newcommand{\erema}{\end{rema}}
\newcommand{\bexer}{\begin{exer}}
\newcommand{\eexer}{\end{exer}}
\newcommand{\bconj}{\begin{conj}}
\newcommand{\econj}{\end{conj}}
\newcommand{\elemm}{\end{lemm}}
\newcommand{\bcoro}{\begin{coro}}
\newcommand{\ecoro}{\end{coro}}
\newcommand{\bexem}{\begin{exem}}
\newcommand{\eexem}{\end{exem}}

\newcommand{\rem}{\noindent{\bf Remarque. }}

\newcommand{\I}{{\mathcal I}}
\newcommand{\T}{{\mathcal T}}

\renewcommand{\L}{{\mathcal L}}

\newcommand{\E}{{\mathcal E}}

\renewcommand{\O}{{\mathcal O}}
\newcommand{\C}{{\mathcal C}}

\newcommand{\F}{{\mathcal F}}

\usepackage{amssymb}
\usepackage{amsmath}

\newcommand{\maths}[1]{{\mathbb #1}}  

\newcommand{\RR}{\maths{R}}
\newcommand{\NN}{\maths{N}}
\newcommand{\CC}{\maths{C}}
\newcommand{\QQ}{\maths{Q}}
\newcommand{\BB}{\maths{B}}
\renewcommand{\SS}{\maths{S}}

\newcommand{\HH}{\maths{H}}

\newcommand{\FF}{\maths{F}}

\newcommand{\ZZ}{\maths{Z}}
\renewcommand{\AA}{\maths{A}}
\newcommand{\PP}{\maths{P}}

\newcommand{\LL}{\maths{L}}
\newcommand{\TT}{\maths{T}}

\newcommand{\SSS}{{\mathfrak S}}

\newcommand{\ra}{\rightarrow}

\newcommand{\ov}[1]{{\overline{#1}}} 
\newcommand{\wt}[1]{{\widetilde{#1}}}

\newcommand{\ga}{\gamma}
\newcommand{\Ga}{\Gamma}

\newcommand{\PSL}[2]{{\operatorname{PSL}_{#1}(#2)}}
\newcommand{\SL}[2]{\operatorname{SL}_{#1}(#2)}
\newcommand{\GL}[2]{\operatorname{GL}_{#1}(#2)}

\newcommand{\Modg}{{\operatorname{Mod}(\Sigma_g)}}
\newcommand{\Modgpm}{{\operatorname{Mod}_\pm(\Sigma_g)}}
\newcommand{\Cg}{{\C(\Sigma_g)}}
\newcommand{\Outn}{{\operatorname{Out}(\FF_n)}}
\newcommand{\Autn}{{\operatorname{Aut}(\FF_n)}}
\newcommand{\Homg}{{\operatorname{Homeo}(\Sigma_g)}}
\newcommand{\Homo}{{\operatorname{Homeo}_0(\Sigma_g)}}
\newcommand{\Homp}{{\operatorname{Homeo}_+(\Sigma_g)}}
\newcommand{\Difo}{{\operatorname{Diff}_0(\Sigma_g)}}
\newcommand{\Difp}{{\operatorname{Diff}_+(\Sigma_g)}}
\newcommand{\Teichg}{{\operatorname{Teich}(\Sigma_g)}}
\newcommand{\cvn}{\operatorname{CV}_n}

\renewcommand\mathcal{\mathscr}

\newcommand{\cinf}{{{\rm C}^\infty}}

\newcommand{\id}{\operatorname{id}}

\renewcommand{\labelenumi}{(\arabic{enumi})}

\newcommand{\Isom}{{\operatorname{Isom}}}

\renewcommand{\theequation}{-\,\arabic{equation}\,-}

\newcounter{fig}

\usepackage{epsfig}
\usepackage[dvips]{color}

\begin{document}
\maketitle

\noindent{\bf INTRODUCTION}

Soient $m,n,g\in\NN$; consid\'erons les trois beaux groupes suivants :
le groupe sp\'ecial lin\'eaire entier $\SL m\ZZ$, vu comme un arch\'etype
de sous-groupe arithm\'etique non uniforme; le groupe $\Modg$ des
classes d'isotopies d'hom\'eomorphismes pr\'eservant l'orientation
d'une surface compacte connexe orientable $\Sigma_g$ de genre $g$
(appel\'e le {\it groupe modulaire}, ou groupe des diff\'eotopies, et
en anglais ``Mapping class group'', de $\Sigma_g$); et enfin le groupe
$\Outn$, quotient du groupe des automorphismes d'un groupe libre
$\FF_n$ de rang $n$ par son sous-groupe des conjugaisons (appel\'e le
{\it groupe des automorphismes ext\'erieurs} de $\FF_n$). 

Le but de ce rapport est de d\'ecrire quelques-unes des analogies bien
connues et fructueuses entre ces trois groupes, et surtout celles plus
prononc\'ees entre $\Modg$ et $\Outn$, qui ont \'et\'e mises en
\'evidence en particulier par les travaux de K.~Vogtmann (voir les
survols \cite{Vogtmann02,Bestvina02,BriVog06}, dont le dictionnaire de
M.~Bestvina dans \cite{Bestvina02} et \cite[\S 4.7]{BesFei10}).  Ces
analogies sont des moteurs actuels de tr\`es nombreux r\'esultats sur
$\Outn$, dus en particulier  \`a Y.~Algom-Kfir, M.~Bestvina,
M.~Bridson, M.~Clay, B.~Farb, M.~Feighn, V.~Guirardel,
U.~Hamenst\"adt, M.~Handel, I.~Kapovich, G.~Levitt, M.~Lustig,
R.~Martin, L.~Mosher, A.~Pettet, K.~Vogtmann. Ces trois groupes
apparaissent dans de nombreux domaines des math\'ematiques, et une
liste exhaustive de leurs similitudes et diff\'erences actuellement
d\'emontr\'ees ne serait pas raisonnable dans le format
imparti. Apr\`es une premi\`ere partie concernant les motivations,
nous nous concentrerons surtout sur des beaux espaces contractiles sur
lequels ces groupes agissent avec des propri\'et\'es analogues, sur
les propri\'et\'es communes de leurs sous-groupes, et sur les
propri\'et\'es semblables (ou conjecturellement semblables) de leur
g\'eom\'etrie asymptotique. Nous renvoyons par exemple \`a
\cite{Vogtmann06} pour une pr\'esentation de propri\'et\'es
cohomologiques communes de ces trois groupes.

{\small\it Le r\'edacteur remercie K.~Vogtmann et G.~Levitt pour des
  commentaires sur une version pr\'eliminaire de ce rapport, et
  ce dernier pour des conversations \`a l'origine de sa pr\'esentation.}

\section{POURQUOI  \'ETUDIER CES GROUPES ?}
\label{sec:motiv}

Pour tout groupe $G$, nous noterons $\operatorname{Out}(G)$, et
appellerons {\it groupe des automorphismes ext\'erieurs} de $G$, le
groupe quotient du groupe $\operatorname{Aut}(G)$ des automorphismes
de $G$ par son sous-groupe (distingu\'e) $\operatorname{Int}(G)=
\{i_g\colon h\mapsto ghg^{-1}\,,\; g\in G\}$ des {\it automorphismes
  int\'erieurs} de $G$. Puisque le groupe fondamental d'un espace
topologique connexe par arcs point\'e ne d\'epend pas du point base
modulo automorphismes int\'erieurs, un choix de point base sera
implicite, et toute application continue $f$ entre deux espaces
topologiques connexes par arcs induit un morphisme $f_*$ entre leurs
groupes fondamentaux, bien d\'efini modulo conjugaison au but.

Pour tous $m,n,g\in\NN$, nous allons nous int\'eresser au groupe des
automorphismes ext\'erieurs de trois groupes~:

\noindent\begin{minipage}{11.2cm}
  $\bullet$~~un groupe ab\'elien libre $\ZZ^m=\pi_1(\TT^m)$ de rang $m$,
  o\`u \mbox{$\TT^m=\prod_{i=1}^m\SS_1$} est un tore de dimension $m$;
\end{minipage} 
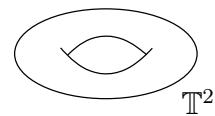
\begin{minipage}{3.7cm}\begin{center}
\begin{picture}(0,0)%
\includegraphics{fig_tore2.pstex}%
\end{picture}%
\setlength{\unitlength}{4144sp}%
\begingroup\makeatletter\ifx\SetFigFontNFSS\undefined%
\gdef\SetFigFontNFSS#1#2#3#4#5{%
  \reset@font\fontsize{#1}{#2pt}%
  \fontfamily{#3}\fontseries{#4}\fontshape{#5}%
  \selectfont}%
\fi\endgroup%
\begin{picture}(1096,704)(623,-400)
\put(1616,-331){\makebox(0,0)[lb]{\smash{{\SetFigFontNFSS{12}{14.4}{\rmdefault}{\mddefault}{\updefault}{\color[rgb]{0,0,0}$\TT^2$}%
}}}}
\end{picture}%
\end{center}
\end{minipage}

\medskip
\noindent \begin{minipage}{11.2cm} $\bullet$~~un groupe libre
  $\FF_n=\pi_1(R_n)$ de rang $n$, o\`u $R_n=\bigvee_{i=1}^n\SS_1$ est un
  bouquet de $n$ cercles point\'es orient\'es; nous noterons $s_1,\dots,
  s_n$ les classes d'homotopies point\'ees de ces cercles, de sorte que
  $\FF_n$ soit l'ensemble des mots r\'eduits en $s_1,\dots, s_n$ et
  leurs inverses, muni de la loi de concat\'enation-r\'eduction;
\end{minipage} 
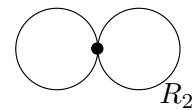
\begin{minipage}{3.7cm}\begin{center}
\begin{picture}(0,0)%
\includegraphics{fig_rose2.pstex}%
\end{picture}%
\setlength{\unitlength}{3729sp}%
\begingroup\makeatletter\ifx\SetFigFontNFSS\undefined%
\gdef\SetFigFontNFSS#1#2#3#4#5{%
  \reset@font\fontsize{#1}{#2pt}%
  \fontfamily{#3}\fontseries{#4}\fontshape{#5}%
  \selectfont}%
\fi\endgroup%
\begin{picture}(1096,715)(623,-409)
\put(1576,-331){\makebox(0,0)[lb]{\smash{{\SetFigFontNFSS{11}{13.2}{\rmdefault}{\mddefault}{\updefault}{\color[rgb]{0,0,0}$R_2$}%
}}}}
\end{picture}%
\end{center}
\end{minipage}

\medskip
\noindent \begin{minipage}{11.2cm} $\bullet$~~un groupe de surface
  $\pi_1(\Sigma_g)$, o\`u $\Sigma_g=\Sigma_{g,0}=\#_{i=1}^g
  \;\SS_1\times\SS_1$ \\ est une somme connexe de $g$ tores de
  dimension $2$, et plus g\'en\'eralement $\Sigma_{g,p}$ est une
  surface lisse compacte connexe orientable de genre $g$ priv\'ee des
  int\'erieurs de $p$ disques ferm\'es plong\'es deux \`a deux
  disjoints, pour tout $p\in\NN$.
\end{minipage} 
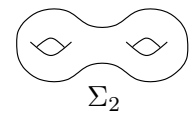
\begin{minipage}{3.7cm}\begin{center}
\begin{picture}(0,0)%
\includegraphics{fig_surfgenre2.pstex}%
\end{picture}%
\setlength{\unitlength}{3729sp}%
\begingroup\makeatletter\ifx\SetFigFontNFSS\undefined%
\gdef\SetFigFontNFSS#1#2#3#4#5{%
  \reset@font\fontsize{#1}{#2pt}%
  \fontfamily{#3}\fontseries{#4}\fontshape{#5}%
  \selectfont}%
\fi\endgroup%
\begin{picture}(1149,734)(619,-1104)
\put(1096,-1026){\makebox(0,0)[lb]{\smash{{\SetFigFontNFSS{11}{13.2}{\rmdefault}{\mddefault}{\updefault}{\color[rgb]{0,0,0}$\Sigma_2$}%
}}}}
\end{picture}%
\end{center}
\end{minipage}

\bigskip
Voici trois raisons d'\'etudier ces groupes et leurs automorphismes.

\smallskip i) Les groupes $\ZZ^m$ et $\FF_n$ v\'erifient des
propri\'et\'es universelles fort utiles : tout groupe (respectivement
groupe ab\'elien) muni d'une partie g\'en\'eratrice ordonn\'ee \`a $n$
(respectivement $m$) \'el\'ements est naturellement un quotient de
$\FF_n$ (respectivement $\ZZ^m$). L'\'etude des groupes
d'automorphismes de tels objets universels dans de telles cat\'egories
est un but en soi.

\smallskip ii) Les surfaces et les hom\'eomorphismes entre surfaces sont
des objets fondamentaux en topologie de petite dimension. 

Avant de justifier cette affirmation, rappelons les propri\'et\'es
suivantes, dues \`a Dehn \cite{Dehn87}, Baer (1924), Nielsen (1927),
Epstein (1966)~: toute \'equivalence d'homotopie de $\Sigma_g$ dans
$\Sigma_g$ est homotope \`a un hom\'eomorphisme, et m\^eme \`a un
diff\'eomorphisme (voir par exemple \cite{Moise77}); deux
hom\'eomorphismes (ou diff\'eomorphismes) de $\Sigma_g$ sont homotopes
si et seulement s'ils sont isotopes (c'est-\`a-dire homotopes \`a
travers des hom\'eomorphismes (ou diff\'eomorphismes)), et si et
seulement si leurs actions sur $\pi_1(\Sigma_g)$ diff\`erent d'un
automorphisme int\'erieur. La premi\`ere propri\'et\'e est
remarquable, et son absence pour le bouquet de cercles $R_n$ (de
nombreux graphes finis ont le m\^eme type d'homotopie que $R_n$, par
exemple) est l'une des raisons pour lesquelles les propri\'et\'es de
$\Outn$ ont des preuves plus techniques que celles pour $\Modg$. Nous
noterons dans la suite $\Homg$ le groupe topologique localement
connexe par arcs des hom\'eomorphismes de~$\Sigma_g$, muni de la
topologie compacte-ouverte, $\Homp$ son sous-groupe topologique
d'indice $2$ des hom\'eomorphismes pr\'eservant l'orientation et
$\Homo$ sa composante neutre, qui est son sous-groupe distingu\'e des
hom\'eomorphismes isotopes \`a l'identit\'e.

{\it D\'ecomposition de Heegaard. } Pour toute vari\'et\'e topologique
$M$ compacte connexe orientable de dimension trois, il existe
$g\in\NN$ et $f:\Sigma_g\ra \Sigma_g$ un hom\'eomorphisme pr\'eservant
l'orientation tel que, en notant $H_g$ un corps \`a anses
de bord $\Sigma_g$, la vari\'et\'e $M$ soit hom\'eomorphe au
recollement $H_g\,\mbox{\tiny $\coprod$}_f H_g$ de deux copies de
$H_g$ le long de $f$. Notons (outre que $g$ n'est pas unique) que l'on
peut prendre pour $H_g$ un voisinage r\'egulier du $1$-squelette d'une
triangulation de $M$, et qu'un tel recollement ne d\'epend que de la
classe d'isotopie de $f$, donc de la classe de $f$ dans $\Homp/\Homo$
(voir par exemple \cite{Scharlemann02}).

{\it Surfaces incompressibles. } Nous renvoyons par exemple \`a
\cite{Jaco80,Kapovich01} pour l'int\'er\^et, de Haken \`a Thurston,
des surfaces compactes connexes plong\'ees dans une vari\'et\'e
compacte connexe de dimension $3$ et induisant une injection sur les
groupes fondamentaux, en particulier des d\'ecoupages le long de
telles surfaces et des recollements par des hom\'eomorphismes (dont
seules les classes d'isotopie importent) le long des surfaces
d\'ecoup\'ees.

{\it Fibr\'es plats en surfaces. } Soient $M$ une vari\'et\'e connexe,
$\wt M\ra M$ un rev\^etement universel de groupe de rev\^etement $\Ga$
et $\rho:\Ga\ra \Homg$ un morphisme. Alors l'espace quotient
$P=\Ga\backslash( \wt M\times \Sigma_g)$ (pour l'action diagonale de
$\Ga$) est une vari\'et\'e qui fibre sur $M$ (via l'application de $P$
dans $\Ga\backslash \wt M$ induite par la premi\`ere projection), de
fibre $\Sigma_g$, qui, \`a hom\'eomorphisme pr\`es, ne d\'epend que de
la classe d'isotopie de $\rho$. Par exemple, si $M=\SS^1$ et
$\Ga=\ZZ$, la g\'eom\'etrie \`a la Thurston (voir par exemple
\cite{Thurston97}) de la vari\'et\'e $P$ de dimension $3$ est dict\'ee
par la classe de $\rho(1)$ dans $\operatorname{Homeo}
(\Sigma_g)/\Homo$. Lorsque $M$ est une surface compacte priv\'ee d'un
nombre fini de points, la compr\'ehension des morphismes de
$\pi_1(M)$ dans $\operatorname{Homeo}(\Sigma_g)/\Homo$ est
importante pour une bonne compr\'ehension des surfaces complexes et des
vari\'et\'es symplectiques de dimension $4$, dont la classification
des fibrations de Lefschetz (voir les travaux de Donaldson, Gompf,
Auroux, et par exemple \cite{KorSti09}).

\smallskip iii) Alors que $\ZZ^m$ est ab\'elien, les groupes $\FF_n$
et $\pi_1(\Sigma_g)$ sont des arch\'etypes de groupes hyperboliques au
sens de Gromov (voir par exemple \cite{Ghys04}). En un certain sens,
ce sont les seuls groupes hyperboliques dont les groupes des
automorphismes ext\'erieurs sont int\'eressants. Bien qu'il reste des
recherches \`a effectuer sur les automorphismes des produits libres,
nous avons en effet les deux r\'esultats suivants. Par des travaux de
Rips et du r\'edacteur (voir par exemple \cite{Paulin97a}), si $\Ga$
est un groupe hyperbolique qui ne se d\'ecompose pas en extension HNN
ou produit amalgam\'e non trivial sur un groupe virtuellement
monog\`ene, alors $\operatorname{Out}(\Ga)$ est fini. Par des travaux
de Sela et Levitt (voir \cite{Levitt04}), si $\Ga$ est un groupe
hyperbolique sans torsion librement ind\'ecomposable, alors il existe
$k,\ell\in\NN$, des surfaces \`a bord (\'eventuellement vide)
compactes connexes $S_1,\dots, S_k$, un sous-groupe~$G$ d'indice fini
de $\operatorname{Out}(\Ga)$ et, en notant
$\operatorname{Homeo}(S_i,\partial S_i)$ le groupe topologique
localement connexe par arcs des hom\'eomorphismes de $S_i$ pr\'eservant
(globalement) chaque composante connexe du bord de $S_i$, une suite
exacte
$$
1\longrightarrow \ZZ^\ell \longrightarrow G \longrightarrow
\prod_{i=1}^k\pi_0\big(\operatorname{Homeo}(S_i,\partial
S_i)\big)\longrightarrow 1\;.
$$

Le but de ce rapport est une \'etude comparative des trois groupes
suivants~:
\begin{itemize}
\item[$\bullet$] le groupe sp\'ecial lin\'eaire entier $\SL m\ZZ$,
  d'indice $2$ dans $\operatorname{Out}(\ZZ^m)=\GL m\ZZ$;
\item[$\bullet$] le groupe modulaire $\Modg$ de $\Sigma_g$, d\'efini
  indiff\'eremment comme
$$
\pi_0(\Homp),\;\;\Homp/\Homo,\;\;\Difp/\Difo\;:
$$
par les r\'esultats initiaux du point (ii) ci-dessus, les morphismes
\'evidents de $\Homp/\Homo$ dans $\pi_0(\Homp)$ et de $\Difp/\Difo$
dans $\Homp/\Homo$ sont des isomorphismes, et ce dernier groupe est
\mbox{d'indice} $2$ dans le groupe $\Modgpm=\Homg/\Homo$, qui est
isomorphe \`a $\operatorname{Out} \big(\pi_1(\Sigma_g)\big)$;
\item[$\bullet$] le groupe $\Outn$ des automorphismes ext\'erieurs de
  $\FF_n$, qui s'identifie avec le groupe des classes d'homotopie des
  \'equivalences d'homotopie du bouquet de cercles~$R_n$ par
  l'application induite sur le groupe fondamental.
\end{itemize}
La litt\'erature sur ces groupes est abondante, et les analogies avec
les groupes modulaires de surfaces sont une justification presque
syst\'ematique des travaux sur $\Outn$.

Les liens directs entre ces trois groupes sont pourtant peu
nombreux. Notons que $\pi_1(\Sigma_{g,p})$ est un groupe libre de rang
$2g+p-1$ si $p>0$, ce qui permet d'\'etudier certains (mais seulement
certains) automorphismes ext\'erieurs de $\FF_{n}$, dits {\it
  g\'eom\'etriques}, tels qu'il existe $g,p\in\NN-\{0\}$ v\'erifiant
$n=2g+p-1$ et un hom\'eomorphisme de $\Sigma_{g,p}$ induisant cet
automorphisme ext\'erieur sur le groupe fondamental, par des
m\'ethodes de topologie et g\'eom\'etrie des surfaces. Les groupes
d'automorphismes de groupes d'Artin \`a angles droits permettent
une certaine interpolation entre les propri\'et\'es des $\GL {m}\ZZ$ 
et des $\Outn$, voir les travaux de Charney-Vogtmann \cite{ChaVog09}.

En petite complexit\'e, nous avons des isomorphismes classiques, le
second \'etant d\^u \`a Nielsen (1918):
$$
\operatorname{Mod}(\Sigma_1)\simeq \SL 2\ZZ\;\;\;{\rm et}\;\;\;
\operatorname{Out}(\FF_2)\simeq \GL 2\ZZ\;.
$$

Le morphisme d'ab\'elianisation $\FF_n\ra\ZZ^n$ induit un morphisme
$\Outn\ra \GL n\ZZ$ qui est surjectif, de noyau de type fini (Magnus
1934), sans torsion (Baumslag-Taylor 1968), qui n'est pas de
pr\'esentation finie si $n=3$ (Krsti\'c-McCool 1997), mais on ne sait
toujours pas s'il est de pr\'esentation finie pour $n\geq 4$.

Le morphisme d'ab\'elianisation $\pi_1(\Sigma_g)\ra \ZZ^{2g}$ (ou
l'action des hom\'eomorphismes en homologie) induit un morphisme
$\Modg\ra \GL {2g}\ZZ$, d'image $\operatorname{Sp}_{2g}(\ZZ)$, de
noyau (appel\'e le {\it groupe de Torelli}) sans torsion, de type fini
si $g\geq 3$ (Johnson 1983), qui n'est pas de type fini si $g=2$
(McCullough-Miller 1986), mais on ne sait toujours pas s'il est de
pr\'esentation finie pour $g\geq 3$.

L'id\'ee g\'en\'erale est que $\Modg$ et $\Outn$ ont assez peu de
morphismes à valeurs dans des groupes lin\'eaires ou venant de groupes
lin\'eaires. Si $n\geq 4$, le groupe $\Outn$ n'admet pas de
repr\'esentation lin\'eaire fid\`ele (Formanek-Procesi 1992), mais on
ne sait pas si $\operatorname{Out}(\FF_3)$ ni si $\Modg$ pour $g\geq
2$ en admet une ou pas. Toutefois, Grunewald et Lubotzky ont construit
dans \cite[Prop.~9.2]{GruLub09} une famille assez riche de
repr\'esentations lin\'eaires (en dimension finie) de sous-groupes
d'indice fini de $\Outn$.  Tout morphisme d'un sous-groupe $\Ga$
d'indice fini de $\SL m\ZZ$ pour $m\geq 3$ (et plus g\'en\'eralement
d'un r\'eseau irr\'eductible $\Ga$ d'un groupe de Lie r\'eel connexe
semi-simple de centre fini, sans facteur compact, de rang r\'eel au
moins $2$), \`a valeur dans $\Modg$ (Farb-Mazur 1998) et dans $\Outn$
(Bridson-Wade 2010), est d'image finie.

\section{ESPACES VIRTUELLEMENT CLASSIFIANTS}
\label{sec:espcalss}

Les groupes $\SL m\ZZ,\Modg,\Outn$ sont de type fini, engendr\'es
respectivement~:

$\bullet$~~par les {\it transformations \'el\'ementaires}
$\left\{ \begin{array}{l}e_i\mapsto e_i+e_k\\e_j\mapsto e_j \;{\rm
      si}\; j\neq i \end{array}\right.$ pour $1\leq i\neq k\leq m$
  o\`u $(e_1,\dots, e_m)$ est la base canonique de $\RR^m$;

\smallskip
\noindent\begin{minipage}{8.5cm}
  ~~$\bullet$~~par les {\it twists de Dehn}, qui sont les classes
  \mbox{d'isotopie} des hom\'eomorphismes de $\Sigma_g$ valant
  l'identit\'e en dehors d'un anneau $\SS_1\times[0,2\pi]$ plong\'e
  dans $\Sigma_g$, et valant
  $(e^{i\theta},t)\mapsto(e^{i(\theta+t)},t)$ sur cet anneau ($2g+1$
  twists de Dehn suffisent \`a engendrer $\Modg$, et ce nombre est
  minimum (Humphries 1979), mais $\Modg$ admet (Wajnryb 1996) une
  paire g\'en\'eratrice);
\end{minipage} 
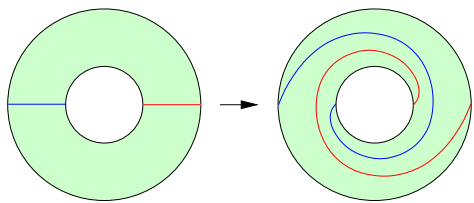
\begin{minipage}{6.4cm}\begin{center}
\begin{picture}(0,0)%
\includegraphics{fig_twistDehn.pstex}%
\end{picture}%
\setlength{\unitlength}{3158sp}%
\begingroup\makeatletter\ifx\SetFigFontNFSS\undefined%
\gdef\SetFigFontNFSS#1#2#3#4#5{%
  \reset@font\fontsize{#1}{#2pt}%
  \fontfamily{#3}\fontseries{#4}\fontshape{#5}%
  \selectfont}%
\fi\endgroup%
\begin{picture}(3616,1514)(2843,-3968)
\end{picture}%
\end{center}
\end{minipage} 

\smallskip 
$\bullet$~~par les {\it transformations de Nielsen}, qui outre les
automorphismes d\'efinis par \mbox{$s_j\mapsto s_{\sigma(j)}$} o\`u
$1\leq j\leq n$, pour les $\sigma\in\SSS_n$, et par
$\left\{\begin{array}{l}s_1\mapsto s_1^{-1}\\ s_j\mapsto s_j \;{\rm
      si}\; 2\leq j\leq n \end{array}\right.$, sont ceux d\'efinis par
$\left\{\begin{array}{l}s_i\mapsto s_is_k\\ s_j\mapsto s_j\;{\rm si}\;
    1\leq j\leq n, j\neq i \end{array}\right.$ pour $1\leq i\neq k\leq
n$ (voir par exemple \cite{LynSch77}), qui redonnent les
transformations \'el\'ementaires par ab\'elianisation.

En fait, ces groupes v\'erifient des propri\'et\'es de finitude bien
plus importantes. Nous allons prendre comme fil directeur de cette
partie \ref{sec:espcalss} le th\'eor\`eme suivant.

\btheo\label{theo:mainespclass} Les groupes $\SL m\ZZ,\Modg,\Outn$
sont de pr\'esentation finie, et de dimension cohomologique virtuelle
respectivement $\frac{m(m-1)}{2}$, $4g-5$, $2n-3$. 
\etheo

Rappelons que la {\it dimension cohomologique virtuelle} d'un groupe
virtuellement sans torsion est la dimension cohomologique commune de
ses sous-groupes d'indice fini sans torsion, qu'elle est minor\'ee par
le rang d'un sous-groupe nilpotent sans torsion, et major\'ee par la
dimension d'un CW-complexe contractile localement fini sur laquelle le
groupe agit proprement avec quotient fini (voir par exemple
\cite{Brown82}). Avec un peu de chance (ce qui sera le cas pour $\SL
m\ZZ $ et $\Outn$ mais pas pour $\Modg$), une strat\'egie de
minoration-majoration permet ainsi le calcul de dimensions
cohomologiques virtuelles.

La premi\`ere assertion est respectivement : due \`a Nielsen (si
$n=3$) et à de S\'eguier \cite{Seguier24} (voir aussi les \'ecrits de
Magnus (1934), Steinberg (1962), Milnor \cite[Coro.~10.3 \&
Theo.~10.5]{Milnor71}); attribu\'ee \`a McCool (1975) bien que le
r\'esultat d\'ecoule du fait que l'espace des modules des courbes
complexes lisses de genre $g$ soit une vari\'et\'e quasi-projective
comme montr\'e par Baily (1960), avec pr\'esentation explicitable par
Hatcher-Thurston (1980) et explicit\'ee par Wajnryb (1983, 1999); et
due \`a Nielsen (1924), avec une pr\'esentation plus longue de
d\'emonstration plus courte par McCool (1974). La seconde assertion est
due respectivement \`a Vorono\"{\i} (1907) (et \`a Borel-Serre (1973)
dans un cadre plus g\'en\'eral), \`a Harer (1986), et \`a
Culler-Vogtmann (1986), ce dernier cas par la construction dont nous
allons parler.

Nous allons construire des espaces topologiques localement compacts
contractiles sur lesquels ces trois groupes agissent proprement, et
indiquer comment les modifier en des CW-complexes localement finis
contractiles sur lesquels ces groupes agissent proprement avec
quotient fini. Ceci implique que ces groupes sont de pr\'esentation
finie (et m\^eme de type VFL, voir \cite[page 200]{Brown82}), et
fournit une majoration de leur dimension cohomologique virtuelle. Nous
donnerons (en indiquant le lien avec les d\'efinitions usuelles) des
d\'efinitions analogues, en tant qu'espaces de modules de structures
marqu\'ees, et en tant qu'espaces d'actions isom\'etriques d'un groupe
fix\'e.

Si $\Ga$ est un groupe et $\E$ un ensemble (de classes d'isom\'etrie
\'equivariante) d'espaces m\'etriques munis d'une action isom\'etrique
de $\Ga$, la {\it topologie de Gromov \'equivariante} (voir
\cite{Paulin88}) est l'unique topologie sur $\E$ dont un syst\`eme
fondamental de voisinages d'un \'el\'ement quelconque $X$ est
$(V_{K,F,\epsilon})_{K,F,\epsilon}$ o\`u $\epsilon>0$, $F$ et $K$ sont
des parties finies de $\Ga$ et~$X$ respectivement, et
$V_{K,F,\epsilon}$ est l'ensemble des $Y\in\E$ tels qu'il existe une
partie finie~$K'$ de $Y$ et une bijection $x\mapsto x'$ de $K$ dans
$K'$ telles que
$$
\forall\;\ga\in F,\;\forall\; x,y\in K,\;\;\;
|\;d(x,\ga y)-d(x',\ga y')\;| <\epsilon\;:
$$
ainsi deux \'el\'ements de $\E$ sont \og proches\fg\ s'ils ont de \og
grandes\fg\ parties finies qui sont \og presque\fg\ isom\'etriques de
mani\`ere \og presque\fg\ \'equivariante par une \og grande\fg\ partie
finie de $\Ga$.

Pour tous $\ga\in \Ga$ et $X\in \E$, l'application
$\operatorname{dis}_{\ga}:X\ra [0,+\infty[$ d\'efinie par
$\operatorname{dis}_{\ga}(x)=d(x,\ga x)$ s'appelle la {\it fonction
  d\'eplacement} de $\ga$ dans $X$, et sa borne inf\'erieure
$\ell_X(\ga)=\inf_{x\in X} d(x,\ga x)\in [0,+\infty[$ s'appelle la
{\it distance de translation} de $\ga$ dans $X$; elle ne d\'epend que
de la classe de conjugaison de $\ga$ dans $\Ga$.  Notons
$\PP\RR_{+}^\Ga$ l'espace topologique quotient de l'espace topologique
produit $[0,+\infty[^{\,_{\operatorname{Int}(\Ga)} \backslash \Ga}$
priv\'e de l'application nulle $0$, par l'action de $\RR^*_+$ par
homoth\'eties. Nous appellerons {\it topologie des distances de
  translation} la topologie sur $\E$ la moins fine rendant continue
l'application de $\E$ dans $\PP\RR_{+}^\Ga$ qui \`a~$X\in \E$ associe
la classe d'homoth\'etie de la fonction distance de translation
$\ell_X$, si $\ell_X\neq 0$ pour tout $X\in \E$. L'action par
pr\'ecomposition de $\operatorname{Aut}(\Ga)$ sur l'ensemble des
fonctions constantes sur les classes de conjugaison de $\Ga$ dans
$[0,+\infty[$ induit une action \`a droite de
$\operatorname{Out}(\Ga)$ sur $\PP\RR_{+}^\Ga$.

\subsection{L'espace sym\'etrique $\E_m$ de $\SL m\RR$}

L'espace sym\'etrique $\E_m$ de $\SL m\RR$ est l'un des espaces
suivants~:
\begin{enumerate}
\item[(i)] l'ensemble $\E'_m$ des classes d'isom\'etrie \'equivariante
  d'actions, propres et libres (donc cocompactes) et de covolume $1$,
  de $\ZZ^m$ sur $\RR^m$ par isom\'etries (donc par translations),
  muni de la topologie de Gromov \'equivariante et de l'action \`a
  droite de $\SL m\ZZ$ par pr\'ecomposition;
\item[(ii)] l'ensemble $\E''_m$ des {\it modules des tores plats
    marqu\'es} de dimension $m$ de volume $1$ -- c'est-\`a-dire
  l'ensemble des couples $(X,f)$, o\`u $X$ est une vari\'et\'e plate
  de volume $1$ et $f:\TT^m\ra X$ un diff\'eomorphisme (appel\'e {\it
    marquage}), quotient\'e par la relation d'\'equivalence $(X,f)\sim
  (X',f')$ si $f'\circ f^{-1}:X\ra X'$ est isotope \`a une
  isom\'e\-trie~--, muni de l'action \`a droite par pr\'ecomposition
  de $\SL m\ZZ$, vu comme le groupe des diff\'eomorphismes lin\'eaires
  pr\'eservant l'orientation de $\TT^m$;
\item[(iii)] l'espace homog\`ene $_{\mbox{$\operatorname{SO}(m)$}}
  \backslash \SL m\RR$, muni de la m\'etrique riemannienne quotient
  d'une m\'etrique riemannienne sur $\SL m\RR$ invariante \`a droite
  par $\SL m\RR$ et \`a gauche par $\operatorname{SO}(m)$, et de
  l'action \`a droite par translations \`a droite de $\SL m\ZZ$ (et
  m\^eme de $\SL m\RR$).
\end{enumerate}

L'application, qui \`a un tore plat marqu\'e $(X,f)$ associe la classe
de conjugaison de la repr\'esentation d'holonomie $\rho:\ZZ^m\ra
\Isom(\RR^m)$ d'une application d\'eveloppante (au sens d'Ehresmann,
voir par exemple \cite[\S 3.4]{Thurston97}) de la m\'etrique plate sur
$\TT^m=\RR^m/\ZZ^m$ image r\'eciproque de la m\'etrique plate de $X$
par le marquage $f$, induit une bijection $\SL m\ZZ$-\'equivariante de
$\E''_m$ dans $\E'_m$. L'action \`a gauche par postcomposition de $\GL
m\RR$ sur l'ensemble des actions propres et libres de $\ZZ^m$ sur
$\RR^m$ par translations, qui est continue et simplement transitive,
induit un hom\'eomorphisme $\SL m\ZZ$-\'equivariant de
$\,_{\mbox{$\operatorname{SO}(m)$}}\backslash\SL m\RR$ dans
$\E'_m$. Nous ne parlerons pas des d\'efinitions de $\E_m$ comme espace
des classes d'homoth\'eties de formes quadratiques d\'efinies positives
sur un espace vectoriel r\'eel de dimension $m$, ni comme espace des
ellipsoïdes de volume $1$ et de centre $0$ dans l'espace euclidien
$\RR^m$.

La vari\'et\'e lisse $\,_{\mbox{$\operatorname{SO}(m)$}} \backslash
\SL m\RR$ est localement compacte et contractile (car diff\'eomorphe
\`a $\RR^{\frac{m(m+1)}2 -1}$ par la d\'ecomposition polaire).
Vorono\"{\i} \cite{Voronoi1907} (ou Ash \cite{Ash77} dans un cadre
plus g\'en\'eral, mais dont la démonstration utilise la connaissance
de la dimension cohomologique virtuelle de $\SL m\ZZ$) a montr\'e que
cette vari\'et\'e se r\'etracte par d\'eformation forte $\SL
m\ZZ$-\'equivariante sur un CW-complexe plong\'e de dimension
$\frac{m(m-1)}2$, invariant et de quotient fini. Or cette dimension
est le rang du sous-groupe unipotent triangulaire sup\'erieur de $\SL
m\ZZ$, qui est nilpotent sans torsion. Le groupe lin\'eaire de type
fini en caract\'eristique nulle $\SL m\ZZ$ \'etant virtuellement sans
torsion par le lemme de Selberg, la strat\'egie susmentionn\'ee permet
donc d'obtenir que la dimension cohomologique virtuelle de $\SL m\ZZ$
est $\frac{m(m-1)}2$.
 
\smallskip
\noindent\begin{minipage}{7cm}
  Lorsque $m=2$, la r\'etraction du demi-plan sup\'erieur $\E_2$ sur
  l'arbre de Bass-Serre de $\SL
  2\ZZ\simeq\ZZ/4\ZZ*_{\ZZ/2\ZZ}\ZZ/6\ZZ$ est facile \`a expliciter
  (voir par exemple \cite[page 53]{Serre83}.
\end{minipage} 
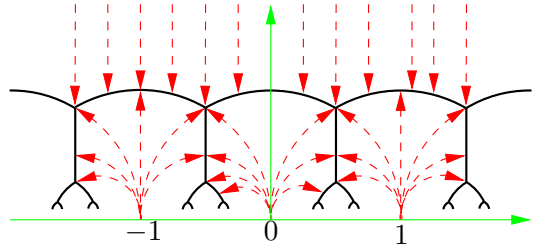
\begin{minipage}{7.9cm}\begin{center}
\begin{picture}(0,0)%
\includegraphics{fig_retractSL2.pstex}%
\end{picture}%
\setlength{\unitlength}{3552sp}%
\begingroup\makeatletter\ifx\SetFigFont\undefined%
\gdef\SetFigFont#1#2#3#4#5{%
  \reset@font\fontsize{#1}{#2pt}%
  \fontfamily{#3}\fontseries{#4}\fontshape{#5}%
  \selectfont}%
\fi\endgroup%
\begin{picture}(3633,1741)(1333,-3440)
\put(2146,-3353){\makebox(0,0)[lb]{\smash{{\SetFigFont{11}{13.2}{\rmdefault}{\mddefault}{\updefault}{\color[rgb]{0,0,0}$-1$}%
}}}}
\put(3106,-3353){\makebox(0,0)[lb]{\smash{{\SetFigFont{11}{13.2}{\rmdefault}{\mddefault}{\updefault}{\color[rgb]{0,0,0}$0$}%
}}}}
\put(4006,-3376){\makebox(0,0)[lb]{\smash{{\SetFigFont{11}{13.2}{\rmdefault}{\mddefault}{\updefault}{\color[rgb]{0,0,0}$1$}%
}}}}
\end{picture}%
\end{center}
\end{minipage} 

\smallskip L'application de $\E'_m$ dans $\PP\RR_{+}^{\ZZ^m}$
d\'efinie par les distances de translation est un hom\'eo\-morphisme
$\SL m\ZZ$-\'equivariant sur son image. Celle-ci est relativement
compacte, et Haettel \cite{Haettel11} a montr\'e que son adh\'erence est
hom\'eomorphe, de mani\`ere $\SL m\ZZ$-\'equivariante, \`a la
compactification de Satake (voir \cite{Satake60}) associ\'ee \`a la
repr\'esentation lin\'eaire tautologique de $\SL m\RR$ sur $\RR^m$
(c'est-\`a-dire \`a l'adh\'erence, dans l'espace projectif des
matrices r\'eelles sym\'etriques $m$-$m$, de l'image de $\SL m\RR$ par
l'application $\ga\mapsto \RR\;^t\ga\ga$), donc est contractile.

\subsection{L'espace de Teichm\"uller $\Teichg$ de $\Sigma_g$}

Fixons une orientation sur $\Sigma_g$ et un plan hyperbolique r\'eel
$\HH^2_\RR$. Si $g\geq 2$, l'espace de Teichm\"uller $\Teichg$ de
$\Sigma_g$ est l'un des espaces suivants~:

\begin{enumerate}
\item[(i)] l'ensemble $\T'$ des classes d'isom\'etrie \'equivariante
  d'actions propres et libres (donc cocompactes) de $\pi_1(\Sigma_g)$
  sur $\HH^2_\RR$, muni de la topologie de Gromov \'equivariante et de
  l'action \`a droite de $\Modgpm$ induite par la pr\'ecomposition par
  l'action sur le groupe fondamental des hom\'eomorphismes de
  $\Sigma_g$;
\item[(ii)] l'ensemble $\T''$ des {\it modules de surfaces
    hyperboliques marqu\'ees} -- c'est-\`a-dire
  l'ensemble des couples $(X,f)$ o\`u $X$ est une surface hyperbolique,
  c'est-\`a-dire munie d'une m\'etrique riemannienne \`a courbure
  sectionnelle constante $-1$
  et $f:\Sigma_g\ra X$ un hom\'eomorphisme
  appel\'e {\it marquage}, quotient\'e par la relation
  d'\'equivalence $(X,f)\sim (X,f')$ si $f'\circ f^{-1}:X\ra X'$ est
  isotope \`a une isom\'e\-trie~--, muni
  de l'action \`a droite de $\Modgpm$ induite par la pr\'ecomposition du
  marquage par les hom\'eomorphismes de $\Sigma_g$;
\item[(iii)] l'espace $\T'''$ des classes d'isotopie de m\'etriques
  hyperboliques (ou de structure analytique complexe induisant
  l'orientation de $\Sigma_g$) sur $\Sigma_g$, muni de la topologie
  quotient de la topologie $\cinf$ sur l'espace des sections lisses du
  fibr\'e des formes bilin\'eaires (ou des endomorphismes) sur
  $\Sigma_g$, et de l'action \`a droite de $\Modgpm$ induite par les
  images r\'eciproques de ces sections par les diff\'eomorphismes de
  $\Sigma_g$.
\end{enumerate}

L'application de $\T''$ dans $\T'''$, qui au module d'une surface
hyperbolique marqu\'ee $(X,f)$ associe la classe d'isotopie de
$F^*\sigma_X$, o\`u $\sigma_X$ est la m\'etrique riemannienne de $X$
et $F:\Sigma_g\ra X$ est un diff\'eomorphisme homotope \`a $f$, est
une bijection $\Modgpm$-\'equivariante. L'application de $\T'''$ dans
$\T'$, qui \`a la classe d'isotopie d'une m\'etrique hyperbolique
$\sigma$ sur $\Sigma_g$ associe la classe de conjugaison de la
repr\'esentation d'holonomie $\rho_\sigma:\pi_1(\Sigma_g)\ra
\Isom_+(\HH^2_\RR)$ d'une application d\'eveloppante de $\sigma$ est
un hom\'eomorphisme $\Modgpm$-\'equivariant (voir par exemple
\cite{Paulin09a}).  Notons que l'espace quotient de $\Teichg$ par
$\Modgpm$ est l'espace des classes d'isom\'etrie de m\'etriques
hyperboliques sur $\Sigma_g$, qui s'identifie, par l'application qui
\`a une structure conforme sur $\Sigma_g$ associe sa m\'etrique
hyperbolique de Poincar\'e, \`a l'espace des modules des structures de
courbe complexe lisse sur $\Sigma_g$. Nous ne parlerons pas ici des
aspects alg\'ebriques et arithm\'etiques de ces espaces de module \`a
la Grothendieck, Baily, Deligne-Mumford, Drinfeld. Par manque de place
(et pour faire plaisir \`a Grothendieck, qui souhaitait \og \'eliminer
un jour \mbox{compl\`etement} l'Analyse de la th\'eorie de l'espace de
Teichm\"uller qui devrait \^etre purement g\'eom\'etrique\fg), nous
n'aborderons pas non plus ses propri\'et\'es d'analyse complexe, voir
par exemple \cite{Gardiner87,Nag88}.

\btheo [Thurston (voir par exemple \cite{FLP})] L'application de $\T'$
dans $\PP\RR_+^{\pi_1(\Sigma_g)}$ d\'efinie par les distances de
translation est un hom\'eomorphisme $\Modgpm$-\'equivariant sur son
image. Celle-ci est relativement compacte, ouverte dans son
adh\'erence, et d'adh\'erence hom\'eomorphe \`a la boule unit\'e
ferm\'ee $\overline{\BB_{6g-6}}$ de $\RR^{6g-6}$, donc contractile.
\etheo

Comme nous avons vu que le groupe de Torelli est sans torsion, et
puisque $\operatorname{Sp}_{2g} (\ZZ)$ est virtuellement sans torsion
par le lemme de Selberg, le groupe $\Modg$ est virtuellement sans
torsion.

La strat\'egie de calcul de la dimension cohomologique virtuelle
employ\'ee pour $\SL m\ZZ$ et $\Outn$ ne fonctionne pas pour $\Modg$~:
la dimension cohomologique virtuelle de $\Modg$ est strictement plus
grande que le rang maximal d'un sous-groupe nilpotent sans torsion de
$\Modg$, et, malgr\'e plusieurs tentatives, il n'a toujours pas
\'et\'e construit de CW-complexe invariant de la bonne dimension
contenu dans $\Teichg$ sur lequel cet espace se r\'etracte par
d\'eformation forte \'equivariante.  Pour montrer que la dimension
cohomologique virtuelle de $\Modg$ est $4g-5$, Harer \cite{Harer86}
utilise alors la suite exacte, pour $x$ un point fix\'e de $\Sigma_g$,
$$
1\longrightarrow\pi_1(\Sigma_g,x)\longrightarrow
\pi_0(\operatorname{Homeo}_+(\Sigma_g,x))\longrightarrow \Modg\;,
$$
ainsi que la strat\'egie de r\'etraction pour l'espace de
Teichm\"uller de $\Sigma_g-\{x\}$ (voir aussi la formulation
hyperbolique par Bowditch-Epstein (1988)), et un lemme de
Bieri-Eckmann (1973) d'additivit\'e par suite exacte de la dimension
cohomologique des groupes sans torsion \`a dualit\'e de Poincar\'e
(notion sur laquelle nous reviendrons plus loin).

\medskip \noindent{\bf Remarque. }  Une autre mani\`ere de construire
la compactification de Thurston de \mbox{l'espace} de Teichm\"uller
est la suivante. Soit $\Ga$ un groupe hyperbolique au sens de
\mbox{Gromov,} non virtuellement monog\`ene. Pour tout automorphisme
$\phi$ de $\Ga$, notons $\widehat\phi$ l'unique hom\'eomorphisme
$\phi$-\'equivariant du bord \`a l'infini $\partial_\infty\Ga$ de
$\Ga$ (pour une partie g\'en\'eratrice de $\Ga$ fix\'ee
indiff\'erente). Notons $\operatorname{Curr}(\Ga)$ l'espace (introduit
par Bonahon (1986) en g\'en\'eralisant une notion de Sullivan) des
{\it courants g\'eod\'esiques} de $\Ga$, c'est-\`a-dire des mesures
bor\'eliennes r\'eguli\`eres positives, invariantes par l'action
diagonale de $\Ga$ et par l'involution $(x,y)\mapsto (y,x)$, sur
l'espace $\partial^2_\infty \Ga$ des couples de points distincts de
$\partial_\infty\Ga$, muni de la topologie vague et de l'action \`a
gauche (commutant avec celle de $\RR^*_+$ par homoth\'eties) de
$\operatorname{Out}(\Ga)$ d\'efinie par $([\phi],\mu)\mapsto
(\widehat\phi\times\widehat\phi)_*\mu$. Notons
$\PP\operatorname{Curr}(\Ga)$ son espace quotient par les
homoth\'eties. Pour tout $\ga\in \Ga-\{e\}$, notons $\mu_\ga$ le
courant g\'eod\'esique, somme des masses de Dirac unit\'es aux couples
de points fixes des conjugu\'es de $\ga$ et de son inverse. L'ensemble
des multiples rationnels des $\mu_\ga$ pour $\ga\in\Ga-\{e\}$ est
dense dans $\operatorname{Curr}(\Ga)$ et $\PP\operatorname{Curr}(\Ga)$
est m\'etrisable compact (voir \cite{Bonahon91} pour ces deux
assertions).

Pour tout espace m\'etrique $X$ g\'eod\'esique (voir définition
ci-dessous) complet localement compact CAT$(-1)$ (voir par exemple
\cite{BriHae98}), muni d'une action isom\'etrique propre de $\Ga$, de
quotient \mbox{compact,} notons $f_X:\partial _\infty
X\ra \partial_\infty \Ga$ l'unique hom\'eomorphisme
$\Ga$-\'equivariant; notons $d_x$ la {\it distance visuelle} sur le
bord \`a l'infini de $X$ vue de $x\in X$, d\'efinie par
$d_x(\xi,\eta)= e^{\frac{1}{2}\lim_{t\ra+\infty}( d(\xi_t,\eta_t)
  -d(x,\xi_t) -d(x,\eta_t))}$, o\`u $t\mapsto \xi_t$ et $t\mapsto
\eta_t$ sont les rayons g\'eod\'esiques de $x$ \`a
\mbox{$\xi\in\partial_\infty X$} et de $x$ \`a $\eta\in\partial_\infty
X$; notons $(\mu_x)_{x\in X}$ une {\it densit\'e de
  Patterson-Sullivan} de $\Ga$ (c'est-\`a-dire une (unique \`a
homoth\'etie pr\`es) famille de mesures positives finies sur le bord
\`a l'infini de~$X$ telles que $\ga_*\mu_x=\mu_{\ga x}$ et
$\frac{d\mu_x}{d\mu_y}(\xi)=e^{-\delta_\Ga\beta_\xi(x,y)}$ pour tous
les $\ga\in\Ga$ et \mbox{$x,y\in X$,}
o\`u~$\delta_\Ga=\limsup_{r\ra+\infty} \frac{1}{n}\ln\operatorname{Card}(\Ga x\cap
B(x,r))$ et \mbox{$\beta_\xi(x,y)= \lim_{t\ra+\infty} d(x,z_t) -
  d(y,z_t)$} o\`u~\mbox{$t\mapsto z_t$} est un rayon g\'eod\'esique
convergeant vers $\xi\in\partial_\infty X$); enfin, notons $m_X$ la
mesure sur $\partial^2 _\infty X$ d\'efinie par $dm_X(\xi,\eta)=
\frac{d\mu_x(\xi)d\mu_x(\eta)} {d_x(\xi,\eta)^{2\delta}}$.

Soit $\E$ un ensemble (de classes d'isom\'etrie \'equivariante)
d'actions isom\'etriques de $\Ga$ sur des espaces métriques, invariant
par homoth\'eties et pr\'ecomposition par $\operatorname{Out}(\Ga)$,
et $\E_0$ son sous-ensemble des actions propres, de quotient compact,
sur des espaces g\'eod\'esiques complets localement compacts
CAT$(-1)$. Notons $\Theta:\E_0\ra \operatorname{Curr}(\Ga)$
l'application qui \`a $X\in \E_0$ associe le courant g\'eod\'esique
$(f_X\times f_X)_*m_X$, qui est \'equivariante (au sens que
$\Theta(xa)=a^{-1}\Theta(x)$ pour tous $x\in \E_0$ et $a\in
\operatorname{Out}(\Ga)$). Un {\it accouplement action-courant} est
une application $\langle\cdot,\cdot\rangle_\E:
\E\times\operatorname{Curr}(\Ga)\ra [0,+\infty[$ continue, homog\`ene
de degr\'e $1$ en chaque variable, telle que $\langle X,\mu_\ga\rangle
= \ell_X(\ga)$ pour tous $X\in \E$ et $\ga\in\Ga$ (elle est unique si
elle existe). Bonahon a montr\'e \cite{Bonahon88}, en prenant
$\Ga=\pi_1(\Sigma_g)$ et $\E=\E_0=\RR_+^*\T'$, outre l'existence d'un
accouplement action-courant, que l'application $\Theta$ induit un
hom\'eomorphisme \'equivariant de $\Teichg$ sur son image dans
$\PP\operatorname{Curr}(\pi_1(\Sigma_g))$, et que l'adh\'erence de son
image est de mani\`ere \'equivariante hom\'eomorphe \`a la
compactification de Thurston de l'espace de Teichm\"uller de
$\Sigma_g$.

\subsection{L'outre-espace $\operatorname{CV}_n$ de 
Culler-Vogtmann}

Avant de d\'efinir l'espace analogue pour $\Outn$ \`a l'espace
sym\'etrique $\E_m$ pour $\SL m\ZZ$ et \`a l'espace de Teichm\"uller
$\Teichg$ pour $\Modg$, rappelons quelques \mbox{d\'efinitions} sur
les arbres r\'eels, introduits par Tits (1977), puis par Morgan et
\mbox{Shalen} (1985) (voir par exemple \cite{Paulin97a} et ses
r\'ef\'erences, dont les survols de Shalen, Morgan, Bestvina).

Un espace m\'etrique $T$ est {\it g\'eod\'esique} si, pour tous $x,y$
dans $T$, il existe une application isom\'etrique d'un intervalle r\'eel
$[a,b]$ dans $T$ envoyant $a$ sur $x$ et $b$ sur $y$. Une {\it
  homoth\'etie} d'un \mbox{espace} m\'etrique $(X,d)$ dans un autre
$(X',d')$ est une application $f:X\ra X'$ telle qu'il existe
$\lambda>0$ tel que $d'(f(x),f(y))=\lambda \,d(x,y)$ pour tous $x,y$
dans $X$. Le groupe $\RR_+^*$ agit par homoth\'eties par
$(\lambda,(X,d))\mapsto (X,\lambda d)$ sur la collection des espaces
m\'etriques. Un {\it graphe m\'etrique} $X$ est un espace m\'etrique
g\'eod\'esique hom\'eomorphe \`a la r\'ealisation g\'eom\'etrique d'un
$1$-complexe simplicial connexe sans sommet terminal. Un {\it sommet}
d'un graphe m\'etrique est un point s\'eparant chacun de ses
voisinages en au moins trois composantes connexes; une {\it ar\^ete}
est l'adh\'erence d'une composante connexe du compl\'ementaire du lieu
des sommets. Un {\it arbre m\'etrique} est un graphe m\'etrique
simplement connexe.

Un {\it arbre r\'eel} est un espace m\'etrique g\'eod\'esique uniquement
connexe par arcs.  Une action isom\'etrique d'un groupe $\Ga$ sur un
arbre r\'eel $T$ est~:

$\bullet$~~ {\it minimale} si elle n'a pas de sous-arbre invariant non
vide propre,

$\bullet$~~ {\it petite} si le stabilisateur de tout arc non trivial
est monog\`ene, 

$\bullet$~~ {\it tr\`es petite} si elle est petite, si aucun
stabilisateur non trivial d'un arc non trivial n'est strictement
contenu dans un groupe monog\`ene, et si aucun \'el\'ement non trivial
de $\Ga$ ne fixe de triplet de points de $T$ non contenus dans un arc.

Notons que Skora (1996) a montr\'e que le bord de Thurston de l'espace
de \mbox{Teichm\"uller} est exactement l'ensemble des classes
d'homoth\'etie des fonctions \mbox{distances} de \mbox{translation}
des actions isom\'etriques minimales petites de $\pi_1(\Sigma_g)$ sur
des arbres r\'eels.

\medskip Si $n\geq 2$, l'{\it outre-espace} $\cvn$ de Culler-Vogtmann,
construit dans \cite{CulVog86}, est l'un des espaces suivants~:

\begin{enumerate}
\item[(i)] l'espace $\cvn'$ quotient, par l'action de $\RR^*_+$ par
  homoth\'eties, de l'ensemble $\operatorname{cv}_n$ des classes
  d'isom\'etrie \'equivariante d'actions isom\'etriques, libres et
  minimales (donc cocompactes) de $\FF_n$ sur des arbres m\'etriques
  au sens ci-dessus, muni de la topologie de Gromov \'equivariante et
  de l'action \`a droite de $\Outn$ induite par la pr\'ecomposition
  des actions par $\operatorname{Aut}(\FF_n)$;
\item[(ii)] l'ensemble $\cvn''$ des {\it modules de graphes
    m\'etriques marqu\'es} -- c'est-\`a-dire l'ensemble des couples
  $(X,f)$ o\`u $X$ est un graphe m\'etrique compact et $f:R_n\ra X$
  une \'equivalence d'homotopie (appel\'ee {\it marquage}),
  quotient\'e par la relation d'\'equivalence $(X,f)\sim (X',f')$ si
  $f'\circ f^{-1}:X\ra X'$ est homotope \`a une homoth\'etie~--, muni
  de l'action \`a droite de $\Outn$ induite par la pr\'ecomposition du
  marquage par les \'equivalences d'homotopie de $R_n$.
\end{enumerate}
L'application qui au module d'un graphe m\'etrique marqu\'e $(X,f)$
associe la classe de l'action de $\FF_n$ sur l'espace total $\wt X$
d'un rev\^etement universel de $X$ (qui est un arbre m\'etrique, en
demandant que ce rev\^etement soit localement isom\'etrique), obtenue
en pr\'ecomposant par $f_*:\FF_n=\pi_1(R_n)\ra \pi_1(X)$ l'action du
groupe fondamental de $X$ sur $\wt X$, est une bijection
\'equivariante de $\cvn''$ dans $\cvn'$. Nous ne parlerons pas de la
d\'efinition de Hatcher-Vogtmann de $\cvn$ utilisant des classes
d'isotopie de sph\`eres de dimension $2$ dans une somme connexe de $n$
copies de $\SS_1\times\SS_2$, pourtant tr\`es bien adapt\'ee pour
certaines applications.

\btheo L'application de $\cvn'$ dans $\PP\RR_+^{\FF_n}$ d\'efinie par
les distances de translation est un hom\'eomorphisme sur son
image. Celle-ci est contractile, \mbox{d'adh\'erence} \mbox{compacte} et
contractile, form\'ee exactement des classes d'homoth\'etie des
fonctions \mbox{distances} de translation des \'el\'ements de l'espace
$\overline{\operatorname{cv}_n}$ des classes d'isom\'etrie
\'equivariante d'actions isom\'etriques minimales tr\`es petites de
$\FF_n$ sur les arbres r\'eels.  \etheo

Cet \'enonc\'e est une combinaison de r\'esultats de Culler-Vogtmann
(1986), Culler-Morgan (1987), Paulin (1989), Cohen-Lustig (1995),
Bestvina-Feighn (non publi\'e), Skora et ind\'ependamment Steiner (non
publi\'es). Point clef de l'article original \cite{CulVog86}, la
d\'emonstration actuellement la plus mall\'eable (en ce qui concerne
les actions de groupes sur les arbres) de la contractibilit\'e de
$\cvn$ qui montre celle de son adh\'erence est sans doute celle de
\cite{GuiLev07a}, suivant les id\'ees de Skora (aussi reprises par
Forester et Clay).

{\small Donnons une id\'ee de cette d\'emonstration. Pour tout
  \'el\'ement $T\in\operatorname{cv}_n$, notons $\wt C(T)$ le
  sous-espace des $T'\in\overline{\operatorname{cv}_n}$ tels que $T'$
  soit hom\'eomorphe, de mani\`ere \'equivariante, \`a l'espace
  quotient de $T$ par l'\'ecrasement en un point de chaque composante
  connexe d'une for\^et invariante de $T$, et $C(T)$ l'image dans
  $\cvn$ de $\wt C(T)\cap \operatorname{cv}_n$. En faisant varier de
  mani\`ere \'equivariante les longueurs des ar\^etes de $T$
  (certaines pouvant devenir nulles), $\wt C(T)$ s'identifie ainsi au
  c\^one \'epoint\'e sur un simplexe de dimension au plus $3n-4$, et
  $C(T)$ \`a ce simplexe priv\'e de certaines facettes ferm\'ees. La
  topologie ainsi d\'efinie sur $\wt C(T)$ co\"{\i}ncide clairement
  avec la topologie de Gromov \'equivariante, et $\wt C(T)$ est donc
  contractile (de mani\`ere \'equivariante par les homoth\'eties).
  Par exemple, si $T$ est un rev\^etement universel de l'un des trois
  (\`a hom\'eomorphisme pr\`es) graphes m\'etriques de groupe
  fondamental isomorphe \`a $\FF_2$, alors $C(T)$ est donn\'e par la
  figure ci-dessous~: un triangle priv\'e de deux c\^ot\'es ferm\'es
  pour le graphe des halt\`eres, un triangle priv\'e de ses sommets
  pour le graphe $\Theta$, un intervalle ouvert pour $R_2$.

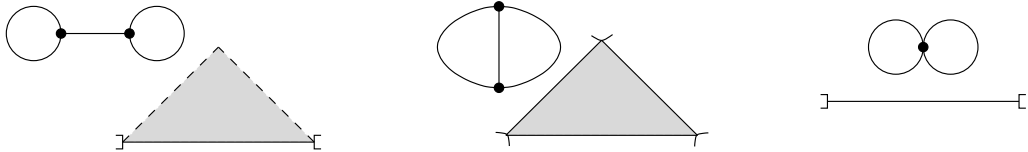
\begin{center}
\begin{picture}(0,0)%
\includegraphics{fig_simplexrang2.pstex}%
\end{picture}%
\setlength{\unitlength}{3729sp}%
\begingroup\makeatletter\ifx\SetFigFontNFSS\undefined%
\gdef\SetFigFontNFSS#1#2#3#4#5{%
  \reset@font\fontsize{#1}{#2pt}%
  \fontfamily{#3}\fontseries{#4}\fontshape{#5}%
  \selectfont}%
\fi\endgroup%
\begin{picture}(6725,996)(-322,-1288)
\end{picture}%

\end{center}

Soit $T_0$ un rev\^etement universel de $R_n$, muni de l'action du
groupe fondamental $\FF_n=\pi_1(R_n)$, et de la distance
g\'eod\'esique rendant chacune de ses ar\^etes isom\'etrique \`a
$[0,1]$.  Munissons $\overline{\operatorname{cv}_n}$ de la topologie
de Gromov \'equivariante et de l'action par homoth\'eties de
$\RR^*_+$. L'id\'ee de Skora est d'expliciter une r\'etraction par
d\'eformation forte $r$ de $\overline{\operatorname{cv}_n}$ sur $\wt
C(T_0)$, pr\'eservant $\operatorname{cv}_n$ et commutant avec l'action
par homoth\'eties de $\RR^*_+$, ce qui conclut.

Si $T$ et $T'$ sont des arbres r\'eels munis d'une action
isom\'etrique d'un groupe $G$, appelons {\it morphisme} de $T$ dans
$T'$ une application \'equivariante telle que tout arc de $T$ puisse
\^etre subdivis\'e en un nombre fini de sous-arcs, en restriction
auxquels $f$ est isom\'etrique. Puisque le graphe d'un morphisme est
un espace m\'etrique muni d'une action isom\'etrique de $G$, nous
munirons de la topologie de Gromov \'equivariante tout ensemble de
morphismes entre \'el\'ements de $\overline{\operatorname{cv}_n}$.
Construisons tout d'abord une application continue $T\mapsto f_T$ de
$\overline{\operatorname{cv}_n}$ dans l'espace des morphismes entre un
\'el\'ement de $\wt C(T_0)$ et un \'el\'ement de
$\overline{\operatorname{cv}_n}$.

Si $\ga$ est une isom\'etrie d'un arbre r\'eel $T$, notons
$A_T(\ga)=\{x\in T\;:\;d(x,\ga x)=\ell_T(\ga)\}$~: alors $A_T(\ga)$
est l'ensemble des points fixes de $\ga$ si celui-ci est non vide, et
sinon $A_T(\ga)$ est isom\'etrique \`a $\RR$ et $\ga$ translate sur
$A_T(\ga)$ de la distance $\ell_T(\ga)>0$.  Rappelons que $s_1$ et
$s_2$ sont les deux premiers \'el\'ements de la famille
g\'en\'eratrice libre fix\'ee de $\FF_n$. Pour tout
$T\in\overline{\operatorname{cv}_n}$, notons $p_T$ le milieu de l'arc
connectant entre $A_{T}(s_1)$ et $A_{T}(s_2)$ s'ils sont disjoints, et
le milieu de leur intersection (qui est un arc de longueur au plus
$\ell_T(s_1)+2\ell_T(s_2)$ car $[s_1,s_2]$ et $[s_1,{s_2}^2]$
engendrent un groupe libre de rang $2$) sinon. Soit $p_*$ le sommet de
$T_0$ intersection de $A_{T_0}(s_1)$ et $A_{T_0}(s_2)$. Notons
$T_0(T)$ l'unique \'el\'ement de $\wt C(T_0)$ tel que l'ar\^ete entre
$p_*$ et $s_ip_*$ soit de longueur $d_{T}(p_T,s_ip_T)$ (non nulle si
$T\in \operatorname{cv}_n$). D\'efinissons alors $f_T:T_0(T)\ra T$
comme l'unique morphisme tel que $f_T(p_*)=p_T$. Il d\'epend
continuement de $T$, est \'equivariant pour l'action par homoth\'eties
de $\RR^*_+$, et vaut l'identit\'e sur $\wt C(T_0)$.

Maintenant, pour construire la r\'etraction $r:[0,1]\times
\overline{\operatorname{cv}_n} \ra \overline{\operatorname{cv}_n}$
voulue, pour tout $t\in[0,1]$ et $T$ dans
$\overline{\operatorname{cv}_n}$, notons $\sim_{t,T}$ la relation
d'\'equivalence sur $T_0(T)$ d\'efinie par $x\sim_{t,T} y$ si et
seulement si $f_T(x)=f_T(y)$ et $d_{T}(f_T(x), f_T(z))\leq -\log (1- t)$
pour tout $z$ dans l'arc entre $x$ et $y$. Notons alors $r(t,T)$
l'ensemble quotient $T_0(T)/\!\sim_{t,T}$ muni de l'action quotient de
$\FF_n$ et de la plus grande distance rendant la projection canonique
$1$-lipschitzienne. L'application $r$ convient (voir par exemple
\cite{GuiLev07b}). }

\medskip Revenons maintenant \`a la d\'emonstration du th\'eor\`eme
\ref{theo:mainespclass} dans le cas de $\Outn$. Comme le noyau du
morphisme canonique $\Outn\ra \GL n\ZZ$ est sans torsion, et puisque
$\GL n\ZZ$ est virtuellement sans torsion par le lemme de Selberg, le
groupe $\Outn$ est virtuellement sans torsion, et sa dimension
cohomologique virtuelle est donc bien d\'efinie.

La topologie de Gromov \'equivariante sur $\cvn$ co\"{\i}ncide (voir
par exemple \cite{GuiLev07a}) avec la topologie faible d\'efinie par
la famille des sous-espaces $C(T)$ pour $[T]\in\cvn$, et $\cvn$ se
r\'etracte de mani\`ere \'equivariante par d\'eformation forte sur la
r\'ealisation g\'eom\'etrique de l'ensemble des $C(T)$ partiellement
ordonn\'e par l'inclusion (voir \cite{CulVog86}). De plus, le quotient
de ce complexe simplicial (qui est donc contractile) par $\Outn$ n'a
qu'un nombre fini de cellules, et il est de dimension $2n-3$, par un
calcul de la caract\'eristique d'Euler d'un graphe m\'etrique
trivalent. Ceci montre que $\Outn$ est de pr\'esentation finie (et
m\^eme de type VFL), et que sa dimension cohomologique virtuelle est
au plus $2n-3$.

Notons que $\Outn$ contient un sous-groupe ab\'elien libre de rang
$2n-3$, l'image dans $\Outn$ du sous-groupe de
$\operatorname{Aut}(\FF_n)$ engendr\'e par les $2n-2$ automorphismes
d\'efinis par $\left\{\begin{array}{l}s_i\mapsto s_is_1\\ s_j\mapsto
    s_j\;{\rm si}\; 1\leq j\neq i\leq n \end{array}\right.$ et
$\left\{\begin{array}{l}s_i\mapsto s_1s_i\\ s_j\mapsto s_j\;{\rm si}\;
    1\leq j\neq i\leq n \end{array}\right.$ pour $2\leq i\leq n$. Ceci
termine le calcul de la dimension cohomologique virtuelle de $\Outn$.

\medskip \rem Alors que la composante neutre du groupe des
isom\'etries de l'espace \mbox{sym\'etrique} $\E_m$, qui est
$\operatorname{PSL}_m(\RR)$, est bien plus grosse que
$\operatorname{PSL}_m(\ZZ)$ si $m\geq 2$, \mbox{Royden} (1971) a
montr\'e que, si $g\geq 2$, le groupe des isom\'etries de $\Teichg$
pour la \mbox{distance} de Teichm\"uller est $\Modgpm$. De m\^eme,
Bridson et Vogtmann (2001) ont montr\'e que l'action de $\Outn$ sur
$\cvn$ induit un isomorphisme de $\Outn$ dans le groupe des
\mbox{automorphismes} simpliciaux du complexe simplicial d\'ecrit
ci-dessus sur lequel se \mbox{r\'etracte} $\cvn$, si $n\geq 3$.

\section{STRUCTURES DES SOUS-GROUPES}
\label{sec:strcutssgroup}

Supposons $m,g,n\geq 2$. Nous prendrons cette fois-ci comme fil
directeur de l'\'etude des sous-groupes de nos trois groupes le
r\'esultat suivant.

\btheo \label{theo:alttits} 
Les groupes $\SL m\ZZ$, $\Modg$ et $\Outn$ v\'erifient
l'alternative de Tits~: tout sous-groupe non virtuellement r\'esoluble
contient un sous-groupe libre de rang~$2$.  
\etheo

Ce r\'esultat a \'et\'e d\'emontr\'e respectivement : par Tits (1972),
voir \cite{Harpe83,Benoist97} pour \mbox{d'excellentes} expositions --
et plus g\'en\'eralement pour tous les groupes lin\'eaires $\Ga$ en
\mbox{caract\'eristique} nulle, et sans hypoth\`ese sur la
caract\'eristique si $\Ga$ est suppos\'e de type fini --; par Ivanov
(1984) et McCarthy (1985); et par Bestvina-Feighn-Handel (2000,
2005). Bien s\^ur, v\'erifier l'alternative de Tits n'est pas
l'apanage de ces trois groupes, voir par exemple les groupes
hyperboliques de Gromov, le groupe des automorphismes polynomiaux du
plan affine complexe (Lamy 2001), ou le groupe de Cr\'emona des
transformations birationnelles du plan projectif complexe (Cantat
2007). Dans le cas de $\Modg$ et $\Outn$, des pr\'ecisions suivantes,
remarquablement semblables, ont \'et\'e apport\'ees respectivement par
Birman-Lubotzky-McCarthy (1983) et Ivanov (1984), et par
Bestvina-Feighn-Handel (2004).

\btheo \label{theo:virtresolmodg}
Tout sous-groupe virtuellement r\'esoluble de $\Modg$ contient un
sous-groupe ab\'elien libre de rang au plus $3g-3$ et d'indice au plus
$3^{4g^2}$.  
\etheo

\btheo\label{theo:virtresoloutn}
Tout sous-groupe virtuellement r\'esoluble de $\Outn$ contient un
sous-groupe ab\'elien libre de rang au plus $2n-3$ et d'indice au plus
$3^{5n^2}$.
\etheo

Nous renvoyons \`a \cite{FeiHan09} pour une description effective de
tous les groupes du th\'eor\`eme \ref{theo:virtresoloutn}.  Comme le
montre la longueur des articles
\cite{BesFeiHan00,BesFeiHan04,BesFeiHan05}, les d\'emonstrations pour
le groupe $\Outn$ sont remarquablement plus techniques que pour les
deux autres groupes, et des ph\'enom\`enes nouveaux apparaissent. Un
d\'enominateur commun est l'\'etude de la dynamique des \'el\'ements
individuels sur les compactifications des espaces construits dans la
partie \ref{sec:espcalss}, ou des variations.  Nous nous concentrerons
sur les \'el\'ements respectivement loxodromiques, pseudo-Anosov,
compl\`etement irr\'eductibles (au sens que nous allons d\'efinir) de
$\SL m\ZZ$, $\Modg$ et $\Outn$, dont la dynamique est la plus simple,
mais qui sont g\'en\'eriques au sens suivant~: une fois fix\'ee une
partie g\'en\'eratrice finie, la proportion d'\'el\'ements de ces
groupes qui sont dans la boule de rayon~$N$ pour la distance des mots
et qui sont respectivement de polyn\^ome caract\'eristique
irr\'eductible sur $\ZZ$, pseudo-Anosov et compl\`etement
irr\'eductibles, tend vers $1$ quand $N\ra+\infty$ (voir
\cite{Rivin08,Kowalski08}). Cette \'etude suffit toutefois pour de
nombreux r\'esultats, par exemple pour la d\'emonstration commune de
la simplicit\'e de l'alg\`ebre stellaire r\'eduite (l'adh\'erence en
norme de l'image de la repr\'esentation r\'eguli\`ere gauche de
l'alg\`ebre complexe du groupe) des quotients par leur centre de ces
trois groupes lorsque $g\geq 3$ et $n\geq 3$, dans \cite{BriHar04}.

\medskip Une technique clef pour construire des groupes libres,
utilis\'ee dans ces trois situations (et d'autres), qui remonte \`a
Klein et Schottky et admet de nombreuses d\'eclinaisons, est la
suivante (voir par exemple \cite{Benoist97}).

\blemm [Lemme du tennis de table] \label{lem:klein} Soit $X$ un espace
topologique compact, muni d'une action par hom\'eomorphismes d'un
groupe $\Ga$. Soit $S$ une partie finie de $\Ga$ telle que $S\cap
S^{-1}=\emptyset$. Supposons que tout $\ga\in S\cup S^{-1}$ admette un
point fixe $x_\ga^+$ attracteur (ceci signifie qu'il existe un
voisinage $U$ de $x_\ga^+$ tel que, pour tout voisinage $V$ de
$x_\ga^+$, nous avons $\ga^n(U)\subset V$ pour tout $n\in\NN$ assez
grand); notons $B^+_\ga=\{x\in X\;:\;\lim_{n\ra+\infty}
\ga^nx=x^+_\ga\}$ le bassin d'attraction (ouvert) de
$x^+_\ga$. Supposons que les $x^+_\ga$ pour $\ga\in S\cup S^{-1}$
soient deux \`a deux distincts et que $x^+_\alpha$ appartienne \`a
$B^+_\beta$ pour tous $\alpha,\beta$ dans $S\cup S^{-1}$ tels que
$\beta\neq \alpha^{-1}$. Alors il existe $p\in\NN-\{0\}$ tel que le
sous-groupe de $\Ga$ engendr\'e par $S^p=\{s^p\;:\;s\in S\}$ soit
libre sur $S^p$.  \elemm

\subsection{\'El\'ements loxodromiques}

Un \'el\'ement de $\SL m\RR$ est {\it loxodromique} si les modules de
ses valeurs propres \mbox{complexes} sont deux \`a deux distincts. Son
action projective sur $\PP_{m-1}(\RR)$ admet alors un point fixe
attractif, correspondant \`a la droite propre de valeur propre de plus
grand module, de bassin d'attraction le compl\'ementaire de
l'hyperplan projectif engendr\'e par les autres droites propres.  La
dynamique projective des autres \'el\'ements de $\SL m\RR$ est
prescrite par sa forme normale de Jordan sur $\CC$.

Une mani\`ere de voir ceci est d'utiliser la fonction de d\'eplacement
$\operatorname{dis}_\alpha$ d'un \'el\'ement~$\alpha$ de $\SL m\RR$
dans l'espace sym\'etrique $\E_m$ muni d'une distance riemannienne
invariante par $\SL m\RR$~: son minimum est atteint si et seulement si
$\alpha$ est semi-simple (diagonalisable sur $\CC$) et lorsqu'il est
atteint, $\alpha$ est loxodromique si et seulement si son ensemble
minimal de translation est exactement un plat maximal de $\E_m$.

Une partie de $\SL m\RR$ est dite {\it unipotente} si ses \'el\'ements
sont trigonalisables de valeurs propres  \'egales \`a $1$ dans une
m\^eme base. Notons le th\'eor\`eme de Kolchin (version pour les groupes
lin\'eaires du th\'eor\`eme d'Engel pour les alg\`ebres de Lie, voir par
exemple 
\cite{Serre92}) qui dit qu'un sous-groupe de $\SL m\RR$ dont tous les
\'el\'ements sont unipotents est unipotent. Une extension de ce
r\'esultat (voir par exemple \cite[Coro.~2.4]{Benoist97}) dit que tout
sous-groupe de type fini de $\SL m\RR$, dont les valeurs propres
complexes de tout \'el\'ement sont des racines de l'unit\'e, est
virtuellement unipotent.

La remarque clef pour montrer le th\'eor\`eme \ref{theo:alttits} pour
$\SL m\ZZ$ est la suivante.

\blemm Un sous-groupe de type fini $\Ga$ de $\SL m\RR$, contenant un
\'el\'ement loxodromique et dont tout sous-groupe d'indice fini agit
de mani\`ere irr\'eductible sur $\RR^m$, contient un sous-groupe libre
de rang $2$.  
\elemm

Il suffit d'appliquer le lemme du tennis de table \`a
$X=\PP_{m-1}(\RR)$ et $S=\{\ga,h\ga h^{-1}\}$ o\`u $\ga\in\Ga$ est
loxodromique et $h\in\Ga$ v\'erifie $hx^+_{\ga^\epsilon}\in
B^+_{\ga^{\epsilon'}}$ pour $\epsilon, \epsilon' =\pm 1$, ce qui est
possible par l'hypoth\`ese d'irr\'eductibilit\'e.

Bien s\^ur, ce lemme ne conclut pas~: $\operatorname{SO}(m)$ par
exemple ne contient pas d'\'el\'ement loxodromique, ni, ce qui suffit
pour la dynamique, d'\'el\'ement ayant, ainsi que son inverse, une
unique valeur propre de module maximal.  Mais si $\Ga$ est un
sous-groupe de type fini de $\SL m\RR$ non virtuellement r\'esoluble,
alors, quitte \`a passer \`a un sous-groupe d'indice fini et en
appliquant l'extension du th\'eor\`eme de Kolchin \`a son sous-groupe
d\'eriv\'e, il contient un \'el\'ement ayant une valeur propre non
racine de l'unit\'e. Celle-ci est de valeur absolue strictement
sup\'erieure \`a $1$ dans au moins un corps local $k$ de
caract\'eristique nulle, contenant les coefficients des \'el\'ements
de $\Ga$, et on remplace $\RR^m$ par $k^m$. Enfin, passer \`a une
puissance ext\'erieure permet d'isoler la plus grande valeur absolue
d'une valeur propre, et un argument de suite de Jordan-H\"older permet
de se d\'ebarrasser de la condition d'irr\'eductibilit\'e.

\subsection{\'El\'ements pseudo-Anosov}

La th\'eorie de Nielsen-Thurston d\'ecrivant les \'el\'ements de
$\Modg$ ayant d\'ej\`a fait l'objet d'un rapport \cite{Poenaru80},
nous ne la rappelons que bri\`evement (voir aussi \cite{FLP}).

Supposons $g\geq 2$. Une {\it multi-courbe} de $\Sigma_g$ est une
union non vide de courbes ferm\'ees simples dans $\Sigma_g$, non
homotopes \`a un point, et deux \`a deux disjointes et non
homotopes. Un \'el\'ement de $\Modg$ est dit

$\bullet$~~ {\it p\'eriodique} s'il est d'ordre fini (ou, de mani\`ere
\'equivalente, s'il admet un repr\'esentant dans
$\operatorname{Homeo}_+(\Sigma_g)$ qui est d'ordre fini);

$\bullet$~~ {\it r\'eductible} s'il admet un repr\'esentant dans
$\operatorname{Homeo}_+(\Sigma_g)$ pr\'eservant une multi-courbe 
de $\Sigma_g$.

$\bullet$~~ {\it pseudo-Anosov} (parce qu'ils g\'en\'eralisent, aux
surfaces $\Sigma_g$ pour $g\geq 2$ qui n'en admettent pas, les
diff\'eomorphismes Anosov, voir par exemple \cite{KatHas95}) s'il existe
un de ses repr\'esentants $h$ dans $\operatorname{Homeo}_+(\Sigma_g)$,
deux feuilletages singuliers transversalement mesur\'es $(\F_+,\mu_+)$
et $(\F_-,\mu_-)$, dits {\it stables et instables}, transverses, et
$\lambda=\lambda(h)\in\;]0,1]$ tels que $h(\F_\pm,\mu_\pm)=
(\F_\pm,\lambda^{\pm 1}\mu_\pm)$.

Les travaux de Nielsen et Thurston montrent que tout \'el\'ement
ap\'eriodique et \mbox{irr\'eductible} de $\Modg$ est pseudo-Anosov,
et d\'ecrivent les \'el\'ements quelconques par une forme normale de
Thurston (voir par exemple \cite{Ivanov01}), qui joue le r\^ole de la
forme normale de Jordan pour l'\'etude dynamique.

Une mani\`ere de voir ceci en utilisant les fonctions de d\'eplacement
est due \`a Bers \cite{Bers78}~: si $\operatorname{dis}_\alpha$ est la
fonction de d\'eplacement d'un \'el\'ement $\alpha$ de $\Modg$ dans
l'espace de Teichm\"uller muni de la distance de Teichm\"uller, alors
son minimum est atteint et nul si et seulement si $\alpha$ est
p\'eriodique; il est atteint non nul si et seulement si $\alpha$ est
pseudo-Anosov; et il n'est pas atteint si et seulement si $\alpha$ est
r\'eductible.

Le point clef de la preuve du th\'eor\`eme \ref{theo:alttits} pour
$\Modg$, qui sert de remplacement \`a l'argument de changement de
corps de base pour $\SL m\ZZ$, et dont la d\'emonstration repose sur
la dynamique de l'action de $\Modg$ sur la compactification de
Thurston, est le suivant.

\bprop [Ivanov \cite{Ivanov92}] \label{prop:ivanovquoi}
Tout sous-groupe infini de $\Modg$ ne pr\'eservant pas de classe
d'isotopie de multi-courbe contient un \'el\'ement pseudo-Anosov.  
\eprop

Un analogue du th\'eor\`eme de Kolchin (qui d\'ecoule de cette
proposition par r\'ecurrence apr\`es adaptation au cas des surfaces
\`a bord) est qu'un sous-groupe de $\Modg$ dont chaque \'el\'ement est
une composition de puissances de twists de Dehn autour des composantes
connexes d'une multi-courbe est constitu\'e de compositions de
puissances de twists de Dehn autour des composantes connexes d'une
m\^eme multi-courbe.

Le caract\`ere hyperbolique de $\pi_1(\Sigma_g)$ et $\FF_n$, que ne
poss\`ede pas $\ZZ^m$ si $m\geq 2$, se refl\`ete sur certaines
propri\'et\'es de l'espace de Teichm\"uller et de l'outre-espace plus
proches des espaces sym\'etriques de type non compact de rang un que
de rang sup\'erieur (voir en particulier les travaux de Minsky sur
l'espace de Teichm\"uller). En particulier, la dynamique des
\'el\'ements g\'en\'eriques de $\Modg$ et $\Outn$ est plus simple.

Un hom\'eomorphisme $h$ d'un espace topologique s\'epar\'e $X$ admet
une {\it dynamique nord-sud} sur $X$ s'il fixe deux points $x,y$ de
$X$ tels que, pour tout voisinage $U$ de $x$ et tout voisinage $V$ de
$y$, pour tout $n\in\NN$ suffisamment grand, $h^n(X-V)\subset U$ et
$h^{-n}(X-U)\subset V$.

Thurston a montr\'e que tout \'el\'ement pseudo-Anosov a une dynamique
nord-sud sur la compactification de Thurston de l'espace de
Teichm\"uller, de points fixes les (classes \mbox{d'homoth\'etie} des
fonctions distances de translation des) arbres réels duaux (au sens,
\`a \mbox{l'origine,} de \cite{MorSha88}) de ses feuilletages stables
et instables, et que le stabilisateur d'un de ces points fixes est
virtuellement monog\`ene.  Ainsi, tout sous-groupe de $\Modg$,
contenant un \'el\'ement pseudo-Anosov $\alpha$ et non virtuellement
monog\`ene, contient un \'el\'ement peudo-Anosov $\beta$ de paire des
points fixes disjointe de celle de $\alpha$. Il contient donc un
groupe libre de rang $2$, en appliquant le lemme du tennis de table
avec $X$ le bord de \mbox{Thurston} de l'espace de Teichm\"uller et
$S=\{\alpha,\beta\}$. Par r\'ecurrence et adaptation au cas des
surfaces \`a bord, le th\'eor\`eme \ref{theo:alttits} d\'ecoule donc
de la proposition pr\'ec\'edente.

\subsection{\'El\'ements compl\`etement irr\'eductibles}

Deux faits marquants rendent l'\'etude dynamique des \'el\'ements de
$\Outn$ (et par l\`a l'\'etude des sous-groupes de $\Outn$) plus
compliqu\'ee que dans le cas de $\Modg$:

$\bullet$~ du point de vue dynamique topologique, alors que tout
\'el\'ement de $\operatorname{Out}(\pi_1(\Sigma_g))$ est induit par
l'action sur le groupe fondamental d'un hom\'eomorphisme de
$\Sigma_g$, un \'el\'ement de $\Outn$ n'est induit par l'action sur le
groupe fondamental que d'une \'equivalence d'homotopie d'un graphe
fini connexe, qui n'est que rarement homotope \`a un hom\'eomorphisme;

$\bullet$~ du point de vue dynamique mesurable, alors que la
compactification de Thurston de l'espace de Teichm\"uller se plonge de
mani\`ere \'equivariante dans l'espace projectifi\'e des courants
g\'eod\'esiques de $\pi_1(\Sigma_g)$, la compactification de
l'outre-espace de Culler-Vogtmann n'admet pas, si $n\geq 3$, de
plongement topologique \'equivariant dans l'espace projectifi\'e des
courants g\'eod\'esiques de $\FF_n$ (voir \cite{KapLus07}).

Le second point pourrait rendre indispensable des \'etudes
parall\`eles et coupl\'ees de la dynamique sur $\overline{\cvn}$ et
$\PP\operatorname{Curr} (\FF_n)$, commenc\'ee par Kapovich et Lustig,
qui ont montr\'e dans \cite{KapLus09} l'existence d'un accouplement
action-courant pour $\Ga=\FF_n$ et $\E= \overline{\operatorname{cv}_n
}$.

Le premier point a rendu n\'ecessaire l'introduction de bonnes
r\'ealisations topologiques des \'el\'ements de $\Outn$, dont les
premiers mod\`eles, d\'ecrits ci-apr\`es, ont \'et\'e d\'efinis par
Bestvina et Handel dans \cite{BesHan92} pour r\'esoudre la conjecture
de Scott (disant que le rang du sous-groupe (libre) des points fixes
d'un automorphisme de $\FF_n$ est au plus $n$, r\'esultat pr\'ecis\'e
dans \cite{GabJagLevLus98}). Les mod\`eles actuellement les plus
\'elabor\'es (les {\it r\'ealisations ferroviaires relatives
  am\'elior\'ees compl\`etement scind\'ees} de \cite{FeiHan09,
  HanMos09pre}) ont des d\'efinitions de plusieurs pages: s'ils sont
utilis\'es pour des descriptions fines \`a sous-groupe d'indice fini
pr\`es de sous-groupes de $\Outn$, leur utilisation est techniquement
lourde.

Un sous-groupe $H$ de $\Outn$ est dit {\it compl\`etement
  irr\'eductible} si aucun sous-groupe d'indice fini de $H$ ne fixe
une classe de conjugaison de facteur libre (voir la définition au
début de la partie \ref{subsec:bordification}) propre de $\FF_n$. Un
\'el\'ement de $\Outn$ (ou de $\operatorname{Aut}(\FF_n)$) est {\it
  compl\`etement irr\'eductible} si le sous-groupe qu'il engendre (ou
son image dans $\Outn$) l'est. Par exemple, l'automorphisme de
$\FF_2$ d\'efini par $\left\{\begin{array}{l} s_1\mapsto s_1 s_2 s_1\\
    s_2\mapsto s_2s_1\end{array}\right.$ est compl\`etement
irr\'eductible, alors que les automorphismes de $\FF_3$ et $\FF_4$
d\'efinis par $\left\{\begin{array}{l} s_1\mapsto s_1,\; s_2\mapsto
    s_1s_2\\ s_3\mapsto s_1s_3s_2 \end{array}\right.$ et
$\left\{\begin{array}{l} s_1\mapsto s_1 s_2 s_1,\; s_2\mapsto s_2s_1\\
    s_3\mapsto s_3s_1s_4s_3,\; s_4 \mapsto s_4s_3
\end{array}\right.$ ne le sont pas.

\medskip Une {\it r\'ealisation ferroviaire} d'un \'el\'ement
compl\`etement irr\'eductible $\alpha$ de $\Outn$ est la donn\'ee d'un
graphe m\'etrique marqu\'e $(X,f)$ et d'une \'equivalence d'homotopie
$\phi:X\ra X$, telle que

$\bullet$~ il existe $\lambda_\alpha>1$ tel que, pour tout $n\in\NN$,
$\phi^n$ soit localement homoth\'etique de rapport ${\lambda_\alpha}^n$
sur chaque ar\^ete de $X$,

$\bullet$~ $\phi$ {\it repr\'esente} $\alpha$, au sens que l'image
dans $\operatorname{Out}(\FF_n)$ de ${f_*}^{-1}\circ\phi_*\circ
f_*:\pi_1(R_n)\ra \pi_1(R_n)$ est $\alpha$.

Leur existence est montr\'ee dans \cite{BesHan92}. L'id\'ee est de
construire un graphe m\'etrique \mbox{marqu\'e} $(X,f)$, une partie
finie $V$ de $X$ contenant les sommets de $X$, une relation
d'\'equivalence pour chaque point $x$ de $X$, ayant au moins deux
classes, sur l'ensemble des germes d'arcs non param\'etr\'es partant
de $x$, et une \'equivalence d'homotopie $\phi:X\ra X$ repr\'esentant
$\alpha$, telle que $\phi$ envoie $V$ dans $V$, est localement
injective en restriction \`a chaque composante connexe de $X-V$, et
envoie deux germes in\'equivalents sur deux germes in\'equivalents. On
munit alors $e\in E=\pi_0(X-V)$ de la longueur $\ell_e$ o\`u
$(\ell_e)_{e\in E}$ est le vecteur propre \`a coordonn\'ees
strictement positives (de norme $\LL^1$ \'egale \`a $1$),
\mbox{associ\'e} \`a la valeur propre $\lambda_{\alpha}$ (strictement
positive) de plus grand module de la \mbox{matrice} irr\'eductible
$(n_{e,e'})_{e,e'\in E}$, o\`u $n_{e,e'}$ est le nombre de points de
$e$ qui s'envoient sur le \mbox{milieu} de $e'$ par $\phi$ (qui existe
par le th\'eor\`eme de Perron-Frobenius). Voir aussi
\cite{Bestvina10pre} pour une construction reposant sur une
minimisation, analogue \`a celle de Bers \'evoqu\'ee pour $\Modg$, de
la fonction de d\'eplacement $\operatorname{dis}_\alpha$ de $\alpha$
sur l'outre-espace muni de la distance non-sym\'etrique
(c'est-\`a-dire v\'erifiant les axiomes d'une distance sauf l'axiome
de sym\'etrie) $d_{\rm Lip}$ suivante~: pour tous $x,y\in \cvn''$, si
$(X,f)$ et $(X',f')$ sont des repr\'esentants de $x$ et~$y$, dont la
somme des longueurs des ar\^etes est $1$, alors $d_{\rm Lip}$ est la
borne inf\'erieure des constantes de Lipschitz des applications
lipschitziennes $\ell:X\ra X'$ telles que $\ell\circ f$ soit homotope
\`a $f'$. L'application $d_{\rm Lip}$ est analogue \`a la distance
non-sym\'etrique de \mbox{Thurston} sur l'espace de Teichm\"uller
(voir par exemple \cite{PapThe07}), et a \'et\'e \'etudi\'ee en
particulier par Whyte, Francaviglia, Martino, Algom-Kfir, Bestvina.

\medskip Le r\'esultat clef pour comprendre la dynamique des
\'el\'ements compl\`etement irr\'eductibles est le suivant, qui
compl\`ete des r\'esultats de Bestvina-Feighn-Handel, Lustig et
Bestvina, et utilise les r\'ealisations ferroviaires.

\bprop [Levitt-Lustig \cite{LevLus03}, Martin \cite{Martin95}] 
Si $n\geq 2$, un \'el\'ement compl\`etement irr\'eductible de $\Outn$
admet une dynamique nord-sud sur l'adh\'erence $\ov{\cvn}$ de
l'outre-espace et, si $n\geq 3$, sur l'unique minimal de
$\PP\operatorname{Curr}(\FF_n)$.  
\eprop

Notons que si $[T_+]$ et $[T_-]$ sont les points fixes attractifs et
r\'epulsifs de $\alpha$ dans $\ov{\cvn}$ (voir aussi
\cite{Paulin97a,Paulin97b} pour leur construction), alors l'action de
$\alpha $ sur $T_\pm$ est $\lambda_{\alpha^{\pm 1}} T_\pm$, o\`u
$\lambda_{\alpha}$ a \'et\'e d\'efini ci-dessus, et de m\^eme en
\'echangeant $\lambda_{\alpha}$ et $\lambda_{\alpha^{- 1}}$ dans
$\PP\operatorname{Curr}(\FF_n)$.  Mais en g\'en\'eral, nous n'avons
pas $\lambda_{\alpha^{-1}}= (\lambda_{\alpha})^{-1}$, contrairement au
cas des \'el\'ements pseudo-Anosov de $\Modg$. De plus, le
stabilisateur de $[T_+]$ dans $\ov{\cvn}$ est virtuellement
monog\`ene, par \cite[Theo.~2.14]{BesFeiHan97}. Comme pour $\Modg$, le
lemme du tennis de table montre alors le r\'esultat principal de
\cite{BesFeiHan97} suivant.

\bprop Un sous-groupe de $\Outn$ contenant un \'el\'ement
compl\`etement irr\'eductible est virtuellement monog\`ene ou contient
un sous-groupe libre de rang $2$.  \eprop

L'analogue du r\'esultat \ref{prop:ivanovquoi} d'Ivanov reste vrai,
bien que de d\'emonstration beaucoup plus technique, comme
mentionn\'e ci-dessus.

\bprop [Handel-Mosher \cite{HanMos09pre}]
Un sous-groupe compl\`etement irr\'eductible de $\Outn$ contient un
\'el\'ement compl\`etement irr\'eductible.  
\eprop

Mais contrairement au cas de $\Modg$ (o\`u un sous-groupe r\'eductible
admet un sous-groupe d'indice fini pr\'eservant chaque composante
connexe du compl\'ementaire d'une multi-courbe, ce qui permet de faire
des r\'ecurrences sur la valeur absolue de la carac\-t\'eristique
d'Euler, au seul prix d'une extension aux surfaces \`a bord des
r\'esultats), un sous-groupe non compl\`etement irr\'eductible de
$\Outn$, s'il admet un sous-groupe \mbox{d'indice} fini pr\'eservant la
classe de conjugaison d'un facteur libre propre $A$, ne pr\'eserve pas
forc\'ement la classe de conjugaison d'un facteur libre $B$ tel que le
morphisme canonique de $A*B$ dans $\FF_n$ soit un isomorphisme. Ceci
justifie l'introduction (voir \cite{BesHan92,BesFeiHan00,BesFeiHan05})
de bonnes r\'ealisations g\'eom\'etriques des \'el\'ements quelconques
de $\Outn$, analogues \`a la forme normale de Jordan d'un \'el\'ement
de $\SL m\RR$ et \`a la forme normale de Thurston d'un \'el\'ement de
$\Modg$. Nous ne le ferons ci-dessous que dans le cas le plus
\'eloign\'e des r\'ealisations ferroviaires d'\'el\'ements
compl\`etement irr\'eductibles.

Pour toute classe de conjugaison $c$ d'un groupe $\Ga$ muni d'une
partie g\'en\'eratrice finie fix\'ee $S$, notons $|c|$ la longueur
minimale d'une \'ecriture comme mot en $S\cup S^{-1}$ d'un \'el\'ement
de $c$. Une partie $H$ de $\operatorname{Out}(\Ga)$ est dite {\it \`a
  croissance des it\'erations polynomiale} si, pour tout $\alpha\in H$
et pour toute classe de conjugaison $c$ de $\Ga$, il existe $k\in\NN$
tel que $\lim_{n\ra+\infty} \frac{|\alpha^n(c)|}{n^{k+1}}=0$.  Par
exemple, un sous-groupe de type fini de $\operatorname{Out}(\ZZ^m)=\GL
m\ZZ$ est \`a croissance des it\'erations polynomiale si et seulement
s'il est virtuellement unipotent, et un sous-groupe de type fini de
$\Modg$ est \`a croissance des it\'erations polynomiale si et
seulement si, \`a sous-groupe d'indice fini pr\`es, il est form\'e de
compositions de puissances de twists de Dehn autour des composantes
connexes d'une m\^eme multi-courbe. Mais un sous-groupe de type fini
de $\Outn$ peut \`a la fois \^etre \`a croissance des it\'erations
polynomiale et contenir un groupe libre de rang $2$. Une {\it
  r\'ealisation ferroviaire triangulaire} d'un \'el\'ement \`a
croissance des it\'erations polynomiale $\alpha$ de $\Outn$ est la
donn\'ee d'un graphe m\'etrique marqu\'e $(X,f)$, d'une partie finie
$V$ de $X$ contenant les sommets de $X$, d'une num\'erotation
$e_1,\dots, e_p$ des adh\'erences des composantes connexes de $X-V$ et
d'une \'equivalence d'homotopie $\phi:X\ra X$ repr\'esentant $f$,
envoyant $V$ dans $V$ et localement injective sur chaque composante
connexe de $X-V$, telle que le chemin $f(e_i)$ soit une
concat\'enation de chemins $u_ie_iv_i$ ou $u_i\overline{e_i}v_i$, o\`u
$u_i,v_i$ sont des lacets dans $\bigcup_{j=1}^{i-1} e_i$.

Bestvina-Feighn-Handel montrent d'une part que l'alternative de Tits
est vraie pour tous les sous-groupes de $\Outn$ si et seulement si
elle est vraie pour les sous-groupes \`a croissance des it\'erations
polynomiale, et d'autre part un analogue du th\'eor\`eme de Kolchin
pour $\SL m\ZZ$ ci-dessus mentionn\'e, disant qu'\`a sous-groupe
d'indice fini pr\`es, tout \'el\'ement d'un sous-groupe de type fini
de $\Outn$ \`a croissance des it\'erations polynomiale admet une
r\'ealisation ferroviaire triangulaire de m\^eme graphe m\'etrique
marqu\'e sous-jacent. Avec une borne sur la complexit\'e de celui-ci,
l'alternative de Tits pour les sous-groupes \`a croissance des
it\'erations polynomiale s'en d\'eduit assez facilement.

Nous renvoyons \`a l'appendice de \cite{BesFeiHan97} pour une
pr\'esentation plus compl\`ete de la structure de la d\'emonstration
dans \cite{BesFeiHan00,BesFeiHan05} du th\'eor\`eme \ref{theo:alttits}
pour $\Outn$.

\section{G\'EOM\'ETRIE ASYMPTOTIQUE}
\label{sec:geomasymp}

Par ce titre vague, nous entendons trois types de r\'esultats sur nos
groupes, concernant la structure \`a l'infini des quotients par ces
groupes des espaces construits dans la partie \ref{sec:espcalss}, la
structure de leurs c\^ones asymptotiques (lorsque ces groupes ou ces
espaces sont regard\'es d'infiniment loin), et certaines propri\'et\'es
de rigidit\'e \`a grande \'echelle. Supposons $m,g,n\geq 2$.

\subsection{Bordification}
\label{subsec:bordification}

Nous appellerons {\it facteur libre} (respectivement {\it facteur
  direct}) d'un groupe $G$ tout sous-groupe propre $A$ de $G$ tel qu'il
existe un sous-groupe propre $B$ de $G$ tel que le morphisme canonique du
produit libre $A*B$ (respectivement du produit direct $A\times B$)
dans $G$ soit un isomorphisme. Nous fixons un rev\^etement universel
$\wt \pi:\wt \Sigma_g\ra \Sigma_g$ de groupe de rev\^etement
$\pi_1(\Sigma_g)$.

\`A l'infini de $\E_m$, $\Teichg$, $\cvn$, vivent des sous-complexes
remarquables invariants par $\SL m\ZZ$, $\Modg$, $\Outn$
respectivement.

$\bullet$~ L'{\it immeuble de Tits sph\'erique} de
$\operatorname{SL}_m$ sur $\QQ$, que nous noterons
$\I_\QQ(\operatorname{SL}_m)$, peut \^etre d\'efini comme la
r\'ealisation g\'eom\'etrique sph\'erique (o\`u l'on remplace chaque
simplexe standard maximal par un simplexe de Coxeter dans la sph\`ere
$\SS_{m-1}$ de type $A_{m-1}$) de l'ensemble partiellement ordonn\'e par
l'inclusion des facteurs directs de $\ZZ^m$, muni de l'action \`a
gauche de $\GL m\ZZ$ induite par l'action lin\'eaire sur les parties
de $\ZZ^m$. L'application, qui \`a un facteur direct de $\ZZ^m$ associe
le stabilisateur dans $\SL m\QQ$ du $\QQ$-sous-espace vectoriel de
$\QQ^m$ qu'il engendre, induit un hom\'eomorphisme simplicial de
$\I_\QQ(\operatorname{SL}_m)$ dans la d\'efinition usuelle de celui-ci
(la r\'ealisation g\'eom\'etrique sph\'erique du complexe simplicial
d'ensemble des sommets l'ensemble $I$ des sous-groupes paraboliques
propres maximaux d\'efinis sur $\QQ$ de $\operatorname{SL}_m$, et de
simplexes les parties $J$ de $I$ telles que $\cap J$ soit
parabolique). En particulier, $\I_\QQ(\operatorname{SL}_m)$ est un
sous-complexe (non localement fini) de l'immeuble de Tits sph\'erique
de $\operatorname{SL}_m$ sur $\RR$ (d\'efini de m\^eme en
rempla\c{c}ant $\QQ$ par $\RR$). Celui-ci admet une (unique) bijection
continue $\operatorname{SL}_m(\RR)$-\'equivariante dans le {\it bord
  g\'eom\'etrique} (l'espace quotient de l'espace des rayons
g\'eod\'esiques par la relation \og \^etre \`a distance de Hausdorff
finie\fg) de la vari\'et\'e riemannienne compl\`ete simplement connexe
\`a courbure sectionnelle n\'egative ou nulle $\E_m=
\,_{\mbox{$\operatorname{SO}(n)$}} \backslash \SL m\RR$.

$\bullet$~ Le {\it complexe des courbes} $\Cg$ de $\Sigma_g$, d\'efini
par Harvey \cite{Harvey79}, est la r\'ealisation g\'eom\'etrique
standard du complexe simplicial (non localement fini) dont les
simplexes sont les ensembles des classes d'isotopie des composantes
connexes d'une multi-courbe, muni de l'action \`a gauche de $\Modg$
induite par l'action de $\Homp$ sur les multicourbes. En
parti\-culier, ses sommets sont les classes d'isotopie des courbes
ferm\'ees simples non homotopes \`a z\'ero de $\Sigma_g$.  Il
s'identifie \`a une partie du bord de Thurston de $\Teichg$, par
l'injection continue qui, au point de coordonn\'ees barycentriques
$(t_1,\dots, t_p)$ du simplexe d\'efini par une multi-courbe $c$ de
composantes connexes $\ga_1,\dots,\ga_p$, associe (la classe
d'homoth\'etie de la fonction distance de translation de)
l'arbre m\'etrique de sommets les composantes connexes de $\wt
\Sigma_g-\wt\pi^{-1}(c)$, d'ar\^etes de longueur $t_i$ les composantes
connexes de $\wt \Sigma_g-\wt\pi^{-1}(\ga_i)$ pour $1\leq i\leq p$, un
sommet $v$ \'etant extr\'emit\'e d'une ar\^ete $e$ si $v$ est contenu
dans l'adh\'erence de $e$, muni de l'action induite (isom\'etrique,
minimale et petite) de $\pi_1(\Sigma_g)$.

$\bullet$~ Le {\it complexe des facteurs libres} $\F_{\rm Out}(\FF_n)$
de $\Outn$ est la r\'ealisation g\'eom\'etrique standard de l'ensemble
des classes de conjugaison des facteurs libres de $\FF_n$
partiellement ordonn\'e par la relation $x\prec y$ s'il existe des
repr\'esentants $A$ et $B$ de $x$ et $y$ tels que $A$ soit un
sous-groupe propre (ou, de mani\`ere \'equivalente, un facteur libre)
de $B$, muni de l'action \`a gauche de $\Outn$ induite par son action
sur les classes de conjugaison des sous-groupes de $\FF_n$. Le
morphisme d'ab\'elianisation de $\FF_n$ dans $\ZZ^n$ induit un
morphisme simplicial surjectif de $\F_{\rm Out}(\FF_n)$ sur l'immeuble de Tits
sph\'erique de $\operatorname{SL}_n$ sur $\QQ$, \'equivariant pour la
projection canonique de $\Outn$ sur $\GL n\ZZ$. Il ne faut pas 
confondre ce complexe avec le {\it complexe des facteurs libres} $\F_{\rm
  Aut}(\FF_n)$ de $\Autn$ d\'efini par Hatcher-Vogtmann (1998), qui est
la r\'ealisation g\'eom\'etrique de l'ensemble partiellement ordonn\'e
par l'inclusion des facteurs libres de $\FF_n$, muni de l'action \`a
gauche induite par l'action de $\Autn$ sur les parties de
$\FF_n$. L'application qui \`a un facteur libre associe sa classe de
conjugaison induit un morphisme simplicial surjectif, \'equivariant
pour la projection canonique de $\Autn$ sur $\Outn$, du second
complexe sur le premier.

Le r\'esultat suivant est d\^u respectivement \`a Solomon-Tits (1969),
\`a Harer (1986) (voir aussi Ivanov (1987)), et \`a Hatcher-Vogtmann
(1998).

\btheo \label{theo:bouquet} L'immeuble de Tits sph\'erique de
$\operatorname{SL}_m$ sur $\QQ$, le complexe des courbes de $\Sigma_g$
et le complexe des facteurs libres de $\Autn$ ont le type d'homotopie
de bouquets de sph\`eres de dimension respectivement $m-2$, $2g-2$,
$n-2$.  
\etheo

La g\'eom\'etrie de $\I_\QQ(\operatorname{SL}_m)$ est bien
comprise. Cet espace m\'etrique de diam\`etre $\pi$ est un immeuble de
Tits de type le syst\`eme de Coxeter sph\'erique de type $A_{m-1}$ et
en particulier est CAT$(-1)$. Pour $k\in\{-1\}\cup \NN$, un espace
topologique localement compact $X$ est {\it $k$-connexe \`a l'infini}
si pour tout compact $K$ de $X$, il existe un compact $K'$ de $X$ tel
que toute application continue de la sph\`ere de dimension $k$ dans le
compl\'ementaire de $K'$ se prolonge continuement en une application
de la boule de dimension $k+1$ dans le compl\'ementaire de $K$. Un
groupe discret est {\it $k$-connexe \`a l'infini} s'il agit
simplicialement, librement, avec quotient compact, sur un complexe
simplicial $k$-connexe \`a l'infini. Un groupe discret est {\it
  virtuellement un groupe \`a dualit\'e de dimension $k$} (une
terminologie introduite par Bieri-Eckmann (1973)) s'il contient un
sous-groupe d'indice fini $\Ga$ agissant cellulairement, librement,
avec quotient compact, sur un CW-complexe contractile et si le groupe
ab\'elien $H^i(\Ga,\ZZ[\Ga])$ est nul pour $i\neq k$ et libre si
$i=k$. Ceci implique que $k$ est la dimension cohomologique virtuelle
du groupe $\Ga$ et que, pour tout $\Ga$-module $A$, les groupes
ab\'eliens $H^i(\Ga,A)$ et $H_{k-i}(\Ga,H^k(\Ga,\ZZ[\Ga])\otimes A)$
sont isomorphes.  Siegel (1953) (g\'en\'eralis\'e par Borel-Serre
\cite{BorSer73}) montre qu'il existe une topologie localement compacte
contractile $\SL m\ZZ$-invariante sur la r\'eunion de $\E_m$ et d'un
espace topologique ayant le type d'homotopie \'equivariant de
$\I_\QQ(\operatorname{SL}_m)$, dont l'espace quotient par $\SL m\ZZ$
est compact. Ce type de construction est depuis connu sous le nom de
{\it bordification}. Il en d\'ecoule que $\SL m\ZZ$ est
$(m-3)$-connexe \`a l'infini, et Borel-Serre (loc. cit.) montrent de
plus qu'il est virtuellement un groupe \`a dualit\'e de dimension
$m(m-1)/2$.

La g\'eom\'etrie de $\Cg$ commence \`a \^etre bien \'etudi\'ee.
Harvey \cite{Harvey79,Harvey81} a construit une bordification de
$\Teichg$ pour le groupe $\Modg$, en rajoutant un espace ayant le type
d'homotopie de $\Cg$ (cette bordification est hom\'eomorphe, de
mani\`ere $\Modg$-\'equivariante, au sous-espace de $\T''$ des
modules de surfaces hyperboliques marqu\'ees de rayon d'injectivit\'e
au moins $\epsilon$, pour tout $\epsilon$ assez petit). Le th\'eor\`eme
\ref{theo:bouquet} montre donc que $\Modg$ est $(2g-3)$-connexe \`a
l'infini, et Harer \cite{Harer86} montre de plus qu'il est
virtuellement un groupe \`a dualit\'e de dimension $4g-5$. 

Le fait le plus marquant est que $\Cg$ est hyperbolique au sens de
Gromov. Ceci est un r\'esultat remarquable (et aux nombreux usages,
que ce soit pour la r\'esolution de la conjecture de la lamination
terminale que pour les propri\'et\'es de la cohomologie born\'ee de
$\Modg$) de Masur-Minsky \cite{MasMin99} (voir \cite{Bowditch06} pour
une d\'emonstration simplifi\'ee et plus effective).  Klarreich (1999)
(voir aussi \cite{Hamenstadt06}) a montr\'e que le bord de Gromov de
$\Cg$ est hom\'eomorphe, de mani\`ere \'equivariante, \`a un espace
quotient d'un sous-espace du bord de Thurston de l'espace de
Teichm\"uller, plus pr\'ecis\'ement l'espace des (classes
d'homoth\'etie des fonctions distances de translation des) actions
libres \`a orbites denses de $\pi_1(\Sigma_g)$, modulo bijections
\'equivariantes pr\'eservant l'alignement. Autrement dit, le bord de
Gromov de $\Cg$ est hom\'eomorphe, de mani\`ere \'equivariante, \`a
l'espace quotient, modulo oubli de la mesure transverse, de l'espace
projectifi\'e des feuilletages singuliers transversalement mesur\'es sur
$\Sigma_g$ qui sont minimaux.

Par contre, la g\'eom\'etrie de $\F_{\rm Out}(\FF_n)$ est actuellement
tr\`es myst\'erieuse. Bestvina-Feighn \cite{BesFei00} ont construit
une bordification de $\cvn$ pour le groupe $\Outn$, par des
consid\'erations de rayons d'injectivit\'e infiniment petits de
diff\'erents ordres pour les graphes m\'etriques marqu\'es, ou, de
mani\`ere \'equivalente, des consid\'erations d'actions de groupes sur
des $\RR^n$-arbres (pour l'ordre lexicographique sur $\RR^n$, voir par
exemple \cite{Chiswell01}). Ils en d\'eduisent que $\Outn$ est
$(2n-5)$-connexe \`a l'infini, et est virtuellement un groupe \`a
dualit\'e de dimension $2n-3$. Mais la relation entre le complexe
$\F_{\rm Out} (\FF_n)$ et le bord de la bordification de
Bestvina-Feighn de l'outre-espace n'est pas claire: il serait
int\'eressant de savoir s'ils ont m\^eme type
d'homotopie. L'hyperbolicit\'e au sens de Gromov de $\F_{\rm Out}
(\FF_n)$ est conjectur\'ee par les experts, mais pour l'instant seuls
des erzatzs ont \'et\'e construits, par Bestvina-Feighn
\cite{BesFei10}. Une possibilit\'e serait d'essayer de construire, par
des techniques d'actions de groupes sur des arbres r\'eels (ou
$\RR^n$-arbres) ou de r\'ealisations ferroviaires et leurs
laminations, le candidat pour \^etre son bord de Gromov \`a l'infini,
puis d'utiliser des m\'ethodes \`a la Tukia, Pansu, Bowditch pour
montrer qu'il est le bord de $\F_{\rm Out} (\FF_n)$ (ou d'un complexe
analogue, voir par exemple \cite{KapLus09} pour une liste de
candidats).

\subsection{C\^ones asymptotiques}
\label{subsec:conasymp}

Soit $(X,d)$ un espace m\'etrique. Choisissons un ultrafiltre sur
$\NN$, plus fin que le filtre de Fr\'echet des compl\'ementaires des
parties finies (voir \cite[I]{Bourbaki71}), $(*_i)_{i\in \NN}$ une
suite de points bases de $X$, et $(\epsilon_i)_{i\in \NN}$ une suite
dans $]0,+\infty[$ convergeant vers $0$. Pour toute suite r\'eelle
$(t_i)_{i\in\NN}$, nous noterons $\lim_\omega t_i$ la limite dans
$[-\infty,+\infty]$ de l'ultrafiltre image de $\omega$ par la suite
(voir \cite[I]{Bourbaki71}).  Consid\'erons l'ensemble
$$
X_\infty=\{x=(x_i)_{i\in \NN}\in X^\NN\;:\; 
\lim_\omega \;\epsilon_i\,d(x_i,*_i)<+\infty\}\;,
$$
point\'e en $(*_i)_{i\in \NN}$, et muni de la pseudo-distance
$d_\infty(x,y) =\lim_\omega \;\epsilon_i\,d(x_i,y_i)$.  Alors l'espace
m\'etrique point\'e $(X_\omega,d_\omega,*_\omega)$ quotient s\'epar\'e
canonique de $(X_\infty,d_\infty,*_\infty)$ est appel\'e un {\it
  c\^one asymptotique} de $X$.

Tout groupe de type fini est muni d'une distance des mots, qui ne
d\'epend pas, ainsi que son c\^one asymptotique \`a autres choix
constants, du choix d'une partie g\'en\'eratrice finie fix\'ee, \`a
hom\'eomorphisme bilipschitz pr\`es. Par exemple, tout c\^one
asymptotique de $\E_2$ et de $\FF_n$ (et donc de $\SL 2\ZZ$,
$\operatorname{Mod}(\Sigma_1)$, $\operatorname{Out}(\FF_2)$ et plus
g\'en\'eralement de tout groupe hyperbolique au sens de Gromov), est
un arbre r\'eel, dont tout point est un point de coupure global.

Si $X$ est un complexe poly\'edral sph\'erique (ayant un nombre fini
de types \mbox{d'isom\'etrie} de cellules, voir \cite{BriHae98} pour
toute information), nous appellerons {\it c\^one radial} sur $X$
\mbox{l'espace} m\'etrique g\'eod\'esique obtenu en rempla\c{c}ant
chaque cellule sph\'erique contenue dans la sph\`ere $\SS_p$ par le
c\^one euclidien des rayons g\'eod\'esiques de $\RR^{p+1}$ qui passent
par ce \mbox{simplexe.} La bordification de l'espace sym\'etrique
$\E_m$ par l'immeuble de Tits sph\'erique de $\operatorname{SL}_m$ sur
$\QQ$ mentionn\'ee ci-dessus a permis \`a Hattori \cite{Hattori95} de
montrer que, pour tout sous-groupe d'indice fini $\Ga$ de $\SL m\ZZ$,
tout c\^one asymptotique de $\E_m/\Ga$ est isom\'etrique au c\^one
radial sur $\Ga \backslash \I_\QQ(\operatorname{SL}_m)$.  Mais une
description aussi pr\'ecise des c\^ones asymptotiques de
$\Teichg/\Modg$ (pour la distance de Teichm\"uller quotient) ne semble
pas \^etre connue.

Les c\^ones asymptotiques de $\SL m\ZZ$ sont assez bien connus. En
effet, par un th\'eor\`eme de Kleiner-Leeb \cite{KleLee97} et une
description plus pr\'ecise de Leeb et Parreau \cite{Parreau10}, tout
c\^one asymptotique de $\E_m$ est isom\'etrique \`a l'immeuble de Tits
affine de $\operatorname{SL}_m$ sur un corps valu\'e complet de
caract\'eristique nulle (muni d'une valuation r\'eelle surjective)
(d\'ependant de l'ultrafiltre $\omega$ et de la suite de scalaires
$(\epsilon_i)_{i\in\NN}$). Par le th\'eor\`eme de
Lubotzky-Mozes-Raghunathan \cite{LubMozRag00}, tout c\^one
asymptotique de $\SL m\ZZ$ s'envoie donc par un hom\'eomorphisme
bilipschitzien dans un immeuble de Tits affine. En particulier, la
dimension topologique maximale d'un sous-espace localement compact
d'un c\^one asymptotique de $\SL m\ZZ$ est $m-1$. Mais on peut montrer
qu'aucun de ses points n'est un point de coupure local si $m\geq 3$.

Les c\^ones asymptotiques de $\Modg$ commencent \`a \^etre compris, en
particulier suite aux travaux de Behrstock. Confirmant que les groupes
modulaires cumulent des propri\'et\'es de rang un et de rang
sup\'erieur, celui-ci a montr\'e dans \cite{Berhrstock06} que tout
point de tout c\^one asymptotique de $\Modg$ est un point de coupure
global.  Behrstock-Minsky \cite{BerMin08} ont montr\'e (en utilisant
les techniques de Minsky de projections fortement contractantes sur
les g\'eod\'esiques de Teichm\"uller dont les images dans le quotient
par $\Modg$ restent dans un compact, qui avaient d\'ej\`a servi dans
\cite{MasMin99}) que la dimension topologique maximale d'un
sous-espace localement compact d'un c\^one asymptotique de $\Modg$ est
$3g-3$ (le m\^eme nombre qui appara\^{\i}t dans le th\'eor\`eme
\ref{theo:virtresolmodg}).

Mais \`a la connaissance de l'auteur, les probl\`emes analogues pour
$\Outn$ sont encore ouverts, en particulier parce que la structure
m\'etrique \`a grande \'echelle de $\cvn$ est obscure. Voir par
exemple le travail d'Algom-Kfir (2009), utilisant la distance non
sym\'etrique sur $\cvn$ \'evoqu\'ee auparavant, pour des r\'esultats
d'existence de points de coupure global dans tout c\^one asymptotique
de $\cvn$, suivant une m\'ethode de d\'emonstration semblable \`a
celle de \cite{Berhrstock06}.

La pr\'esence de points de coupures globaux dans les c\^ones
asymptotiques ne suffit quand m\^eme pas pour montrer
l'hyperbolicit\'e relative (au sens fort, une notion due \`a Gromov,
et d\'evelopp\'ee par Farb et Bowditch) de nos groupes. En fait,
Karlsson-Noskov (2004), Anderson-Aramayona-Shackleton (2007) et
Behrstock-Dru\c{t}u-Mosher (2009) ont montr\'e que si $m\geq 3$,
$g\geq 2$ et $n\geq 3$, alors $\SL m\ZZ$, $\Modg$ et $\Outn$ ne sont
relativement hyperboliques par rapport \`a aucune famille finie de
sous-groupes de type fini d'indice infini.

\subsection{Rigidit\'e}
\label{subsec:rigidite}

Le r\'esultat suivant est, respectivement, une cons\'equence du
th\'eor\`eme de rigidit\'e de Mostow-Margulis (voir par exemple
\cite{Zimmer84}), d\^u \`a Ivanov (voir par exemple
\cite[8.5.B]{Ivanov01}) et d\^u \`a Farb-Handel \cite{FarHan07}
(am\'eliorant fortement le r\'esultat de Kramtsov (1990) que
$\operatorname{Out}(\operatorname{Out}(\FF_n))=\{1\}$ (montr\'e par
Bridson-Vogtmann (2000) en utilisant l'outre-espace)). Dans ces trois
cas, la d\'emonstration utilise, du moins en esprit, les
propri\'et\'es des sous-groupes ab\'eliens (maximaux) et de la
combinatoire de leurs intersections, que ce soit par la combinatoire
des plats maximaux d\'ecrite par l'immeuble de Tits \`a l'infini de
$\SL m\RR$, par les nombres d'intersection de courbes ferm\'es simples
dans $\Sigma_g$ (autour desquelles ont lieu les twists de Dehn), ainsi
que par l'\'etude de \cite{FeiHan09} des sous-groupes ab\'eliens de
$\Outn$ et le th\'eor\`eme \`a la Kolchin de \cite{BesFeiHan05}
mentionn\'e ci-dessus.

\btheo
Soit $\Ga$ un sous-groupe d'indice fini de $\SL m\ZZ$ pour $m\geq 3$,
$\Modg$ pour $g\geq 2$, ou $\Outn$ pour $n\geq 4$. Alors
$\operatorname{Out}(\Ga)$ est fini.
\etheo

Des m\'ethodes semblables fournissent le r\'esultat suivant. Le {\it
  commensurateur abstrait} $\operatorname{Comm}(\Ga)$ d'un groupe
$\Ga$ est le groupe des classes d'\'equivalence d'isomorphismes entre
sous-groupes d'indice fini de $\Ga$, o\`u deux tels isomorphismes
sont identifi\'es s'ils co\"{\i}ncident sur un sous-groupe
d'indice fini de $\Ga$, muni de la loi $[f':G'\ra H' ] \circ[f:G\ra
H]=[f'\circ f:f^{-1}(H\cap G')\ra f'(G'\cap H)]$. Nous avons un
morphisme de $\Ga$ dans $\operatorname{Comm}(\Ga)$, qui \`a
$\ga\in\Ga$ associe la classe de la conjugaison $i_\ga:\Ga\ra \Ga$. Par
exemple, l'application qui, \`a $\ga\in\SL m\QQ$ associe la classe de
la conjugaison par $\ga$ d\'efinie sur le sous-groupe d'indice fini
$\ga^{-1}\SL m\ZZ\ga\cap\SL m\ZZ$, est un isomorphisme de $\SL m\QQ$
sur $\operatorname{Comm}(\SL m\ZZ)$ pour $m\geq 3$, par le
th\'eor\`eme de densit\'e de Borel (voir par exemple
\cite{Zimmer84}). Ce contraste important avec le r\'esultat suivant
pourrait se comprendre par les faits d\'ej\`a mentionn\'es que le
groupe des isom\'etries de l'espace de Teichm\"uller est r\'eduit \`a
$\Modgpm$ si $g\geq 2$, et le groupe des automorphismes combinatoires
de l'outre-espace est r\'eduit \`a $\Outn$ si $n\geq 3$.

\btheo [Ivanov; Farb-Handel] Le morphisme canonique
$\Modgpm\ra\operatorname{Comm}(\Modg)$ est un isomorphisme pour $g\geq
3$; le morphisme canonique $\Outn\ra\operatorname{Comm}(\Outn)$ est un
isomorphisme pour $n\geq 4$.  
\etheo

Rappelons que deux espaces m\'etriques $X$ et $Y$ sont {\it
  quasi-isom\'etriques} s'il existe $\lambda\geq 1$, $c\geq 0$ et une
application $f:X\ra Y$ telle que $-c+\frac{1}{\lambda} \,d(x,y)\leq
d(f(x),f(y))\leq \lambda \,d(x,y)+c$ pour tous $x,y\in X$, et que
$d(y,f(X))\leq c$ pour tout $y\in Y$. Nous renvoyons \`a \cite{Drutu09}
pour des motivations et une pr\'esentation d\'etaill\'ee du probl\`eme
de savoir quels sont les groupes de type fini quasi-isom\'etriques \`a
un (bel) espace m\'etrique donn\'e, et en particulier de savoir si un
groupe de type fini quasi-isom\'etrique \`a un (beau) groupe de type
fini lui est isomorphe, \`a sous-groupes d'indice fini et quotients de
noyau fini pr\`es. Le r\'esultat suivant est d\^u \`a Eskin
\cite{Eskin98} (voir aussi \cite{Drutu00}) pour $\SL m\ZZ$ et \`a
Hamenst\"adt \cite{Hamenstadt07pre} (voir aussi
\cite{BerKleMinMos09}) pour $\Modg$.  En notant $Z(G)$ le centre d'un
groupe $G$, nous avons $Z(\SL m\ZZ)=\{\pm\id\}$ si $m$ est pair, $Z(\SL
m\ZZ)=\{1\}$ sinon; $Z(\Modgpm)\simeq \ZZ/2\ZZ$ si $g=1,2$ (engendr\'e
par la classe de l'involution hyperelliptique si $g=2$) et
$Z(\Modgpm)=\{1\}$ sinon; $Z(\Outn)\simeq \ZZ/2\ZZ$ si $n=1,2$ et
$Z(\Outn)=\{1\}$ sinon.

\btheo Si $G=\SL m\ZZ$ et $m\geq 3$ ou si $G=\Modgpm$ et $g\geq 2$, alors
tout groupe de type fini quasi-isom\'etrique \`a $G$ admet un
sous-groupe d'indice fini qui s'envoie dans $G/Z(G)$ avec noyau fini
et image d'indice fini.  
\etheo

Mais \`a la connaissance du r\'edacteur, le probl\`eme
analogue pour $\Outn$ est encore ouvert.

\medskip Un probl\`eme de rigidit\'e voisin est le suivant. Gromov
\cite[0.2.C]{Gromov93} a montr\'e que deux groupes de type fini sont
quasi-isom\'etriques si et seulement s'ils admettent deux actions
continues, propres et \`a quotient compact, qui commutent, sur un
m\^eme espace topologique localement compact.  Comme d\'efini par
Gromov, deux groupes d\'enombrables $G$ et $H$ sont dit {\it
  mesurablement \'equivalents} s'ils admettent deux actions mesurables
pr\'eservant la mesure, essentiellement libres, admettant un domaine
fondamental mesurable de mesure finie, et qui commutent, sur un m\^eme
{\it espace de Lebesgue} (espace bor\'elien standard muni d'une mesure
positive $\sigma$-finie). Par exemple, deux {\it r\'eseaux}
(sous-groupes discrets de covolume fini) du groupe unimodulaire $\PSL
m\RR$ sont mesurablement \'equivalents (en consid\'erant les actions
par translations \`a droite et \`a gauche), mais ne sont pas
forc\'ement commensurables. Le premier r\'esultat suivant est d\^u \`a
Furman \cite{Furman99} (et plus g\'en\'eralement pour les r\'eseaux de
groupes de Lie r\'eels connexes quasi-simples de centre fini, sans
facteur compact, de rang r\'eel au moins $2$) et le second \`a Kida
\cite{Kida10}.

\btheo Tout groupe d\'enombrable mesurablement \'equivalent \`a $\SL
m\ZZ$, o\`u $m\geq 3$, admet un morphisme dans $\PSL m\RR$ de noyau fini
et d'image un r\'eseau.  Tout groupe d\'enombrable mesurablement
\'equivalent \`a $\Modg$, o\`u $g\geq 2$, admet un morphisme dans
$\Modgpm/Z(\Modgpm)$ de noyau fini et d'image d'indice fini.  
\etheo

Mais \`a la connaissance du r\'edacteur et de Gaboriau, le probl\`eme
analogue pour $\Outn$ est aussi encore ouvert.

\end{document}